\documentclass[12pt]{amsart}
\usepackage{amssymb}
\usepackage{amsfonts} 
\usepackage{amssymb}
\usepackage{multicol}
\usepackage{latexsym}
\usepackage{amsmath}
\usepackage{color}
\usepackage{fullpage}
\usepackage{hyperref}
\usepackage{graphics}
\usepackage{graphpap}
\usepackage{pict2e}
\usepackage{color}
\usepackage{multimedia}
\usepackage{xmpmulti}
\usepackage{tikz}
\usepackage{calligra}
\usepackage{hyperref}
\usepackage{geometry}
\usepackage{amssymb}
\usepackage{mathtools}
\usepackage{tikz}
\usetikzlibrary{chains}

\tikzset{node distance=2em, ch/.style={circle,draw,on chain,inner sep=2pt},chj/.style={ch,join},every path/.style={shorten >=4pt,shorten <=4pt},line width=1pt,baseline=-1ex}

\newcommand{\e}{\mathfrak{e}}
\newcommand{\f}{\mathfrak{f}}

\newcommand{\g}{\mathfrak{g}}
\newcommand{\h}{\mathfrak{h}}
\renewcommand{\k}{{\mathfrak k}}
\renewcommand{\l}{\mathfrak l}
\newcommand{\p}{\mathfrak p}

\newcommand{\s}{\mathfrak s}

\newcommand{\z}{\mathfrak z}
\renewcommand{\t}{{\mathfrak t}}
\renewcommand{\u}{{\mathfrak u}}

\begin{document}
\newtheorem{thm}{Theorem}
\newtheorem{pro}{Proposition}
\newtheorem{cor}{Corollary}
\newtheorem{lem}{Lemma}

\title{Associated symmetric pair and multiplicities of admissible restriction of Discrete Series}
\author{ Jorge A.  Vargas}
\thanks{Partially supported by FONCYT, CONICET, AgenciaCbaCiencia, SECYTUNC
(Argentina), American Institute of Mathematics (USA),
  University of Paris VII (France),    ICTP
(Italy) }
\date{\today }
\keywords{ Branching laws, admissible restriction, Kostant-Blattner's formulae}
\subjclass[2010]{Primary 22E46; Secondary 17B10}
\address{ FAMAF-CIEM, Ciudad Universitaria, 5000 C\'ordoba, Argentine}
\email{vargas@famaf.unc.edu.ar}


\begin{abstract}

Let $( G,   H)$ be a symmetric  pair for a real semisimple Lie group $G$ and $(G, H_0)$   its associated pair. For each   irreducible square integrable representation $\pi$ of $G$ so that its restriction to $H$ is admissible, we find an irreducible square integrable representation $\pi_0$ of $H_0$ which allows to compute the Harish-Chandra parameter of each irreducible $H-$subrepresentation of $\pi$ as well as its multiplicity. The computation   is based on the spectral analysis of the restriction of  $\pi_0$ to a maximal compact subgroup of $H_0.$
\end{abstract}
\maketitle
\markboth{Jorge Vargas}{Dual multiplicities}
\newtheorem{rem}{Remark}
\section*{Introduction}\label{sec:intro}

Branching laws, that is, to write as explicit as possible the decomposition of a given representation in terms of irreducible objets is of importance in several branches of mathematics as well as physics and chemistry. A particular case is to consider an irreducible unitary representation $\pi$ of a group $G$ and a closed subgroup $H$ of $G$   and we wish to find the irreducible $H-$subrepresentatios and the weakly contained irreducible factors, as well as their respective multiplicities. To solve this problem involves among other branches of mathematics,  algebraic geometry, differential geometry and hard analysis as we can learn from  examples and theorems presented in \cite{kobic},  references therein and further work of T. Kobayashi, N. Wallach as well as other researchers. In the book  \cite{knapp}, or in \cite{kita}, \cite{kobmem}, \cite{sato} as well as in the work of other authors \cite{gw}, we learn that sometimes the problem of writing the branching law for $\pi$ is translated into the problem of computing branching law for another pair of groups  $ L \subset H_0$ and certain irreducible representation of $H_0.$  From this point of view, in  \cite{os}, \cite{Speh}, \cite{gw} is analyzed the restriction of  a family of Zuckerman modules for a real reductive Lie group $G$ and  $H$ the connected component of the fix point group of an involution of $G.$   Henceforth,  $G$ denotes a connected {\it simple} matrix Lie group. We fix a maximal compact subgroup $K$ of $G$ and   a maximal torus $T$ in $K.$  Harish-Chandra showed that $G$ admits square integrable irreducible representations if and only if $T$ is a Cartan subgroup of $G.$  For this note,  we always assume $T$ is a Cartan subgroup of $G.$ Under these hypothesis, Harish-Chandra  showed that  the set of equivalence classes of irreducible square integrable representations is parameterized by a lattice in $i\mathfrak t_\mathbb R^\star.$ In order to state our results we need to explicit this parametrization and set up some notation.   As usual,  the Lie algebra of a Lie group is denoted by the corresponding lower case German letter followed by the subindex $\mathbb R.$  The complexification of  the Lie algebra of a Lie group is  denoted by the corresponding German letter without any subscript.     $V^\star $ denotes the dual space to a vector space $V.$ Let $\theta$ be the Cartan involution which corresponds to the subgroup $K,$ the associated Cartan decomposition is denoted by $\g=\k +\p.$ Let $\Phi(\g,\t) $ denote the root system attached to the Cartan subalgebra $\t.$ Hence, $\Phi(\g,\t)=\Phi_c \cup \Phi_n =\Phi_c(\g, \t) \cup \Phi_n (\g, \t)$ splits up as the union the set of compact roots and the set of noncompact roots. From now on, we fix a system of positive roots $\Delta $ for $\Phi_c.$   For this note, either the highest weight or the infinitesimal character of an irreducible representation of $K$ is  dominant with respect to $\Delta.$ The Killing form gives rise to an inner product $(..,..)$ in $i\t_\mathbb R^\star.$ As usual, let $\rho=\rho_G $ denote half of the sum of the roots for some system of positive roots for $\Phi(\g, \h).$  \textit{A Harish-Chandra parameter} for $G$ is $\lambda \in i\t_\mathbb R^\star$ such that $(\lambda, \alpha)\not= 0 , $ for every $\alpha \in \Phi(\g,\t) ,$    and so that $e^{\lambda + \rho}$ is a character of $T.$ To each Harish-Chandra parameter, Harish-Chandra associates a unique irreducible square integrable representation $\pi_\lambda^G$ of $G.$ Moreover, he showed the map $\lambda \rightarrow \pi_\lambda^G$  is a bijection from the set of Harish-Chandra parameters dominant with respect to $\Delta$  onto the set of equivalence classes of irreducible square integrable representations for $G.$

 Each Harish-Chandra parameter $\lambda$ gives rise to a system of positive roots
$$\Psi_\lambda =\Psi_{G,\lambda} =\{  \alpha \in \Phi(\mathfrak g, \mathfrak t) : (\lambda, \alpha) >0 \}.$$ From now on, we assume that Harish-Chandra parameter for $G$ are dominant with respect to  $\Delta.$ Whence,  $\Delta \subset \Psi_\lambda.$

  To follow, we fix  a nontrivial  involutive  automorphism $\sigma$ of $G.$ After we replace $\sigma$ by some conjugate we may and will assume $\sigma$ commutes with $\theta.$ Thus, $\sigma \theta$ is another involution of $G$ which commutes with $\theta.$  We write $$\mathfrak h:=\{ X \in \g : \sigma (X)=X\}, \,\,\,\,\mathfrak q :=\{ X \in \g : \sigma (X)=-X \}, \h_0 =\{X \in \g : \sigma \theta (X)=X \}$$ and $H$ (resp. $H_0$)  the identity  connected component of the set of fix points of $\sigma$ (resp. $\sigma \theta)$ in $G.$ Then $(G,H)$ is a symmetric pair as well as $(G,H_0),$ the later pair is called the associated pair to the former pair. Certainly, $H, H_0$ are closed reductive subgroups of $G$ and the following decompositions hold: $$  \g=\h \oplus \mathfrak q, \,   \h =\h \cap \k + \h \cap \mathfrak p,\, \mathfrak q = \mathfrak q \cap \k + \mathfrak q \cap \p, \,  \h_0 = \h \cap \k + \p \cap \mathfrak q. $$ We notice that $H$ and $H_0$ share the maximal compact subgroup $L:= H\cap K=H_0\cap K.$ We may and will assume $T$ is invariant under the involution $\sigma,$ hence $U:=T \cap H$ is maximal torus of $L.$ Since c.f. (0.2)  we  dealt with irreducible square integrable representations of $G$ with an admissible restriction to $H,
  $ without lost of generality we may  assume
 \begin{center} $U$ is a Cartan subgroup of $H$
 \end{center}
 and hence $U$ is a compact Cartan subgroup of $H_0.$ We verify $\Psi_\lambda$ is invariant under the involution $\sigma$. Therefore, there exists a vector in $i\u_\mathbb R $ which is regular and dominant for $\Psi_\lambda,$ in turn, this vector determines  system of positive roots
  \begin{center}
   $\Delta_0, $  $\Psi_{H,\lambda}, \Psi_{H_0,\lambda} $ for respectively $\Phi(\l, \u),   \Phi(\h, \u), \Phi(\h_0, \u) $
  \end{center}
  \noindent
    such that  for  $\alpha \in \Psi_\lambda, $  $\alpha$ restricted to $\u$  takes on  positive values in the Weyl chamber for $\Delta_0.$ Whenever, $U=T$ we have $ \Delta_0 \subset \Delta, \Psi_{H, \lambda} =\Phi(\h, \t)\cap \Psi_\lambda, \Psi_{H_0, \lambda} =\Phi(\h_0, \t)\cap \Psi_\lambda.$

\smallskip

\noindent
  (0.1) For this note, Harish-Chandra parameters for $H $ as well as for $ H_0$ and $L$ are dominant with respect to $\Delta_0.$
\smallskip

\noindent
(0.2) To continue, we assume $\pi_\lambda^G$ restricted to $H$ is an admissible representation.

 \smallskip

 This hypothesis yields   there exits a family of irreducible square integrable representations, $\pi_\mu^H,$ of $H,$ for which the set of corresponding Harish-Chandra parameters is denoted by $$ Spec_H(\pi_\lambda^G)= \{ \mu \in i\u_\mathbb R^\star : \pi_\mu^H \hookrightarrow res_H(\pi_\lambda^G) \}, $$ together with a sequence of positive integers $m^H(\lambda, \mu),  \mu \in Spec_H(\pi_\lambda^G),  $ so that $$ res_H(\pi_\lambda^G) = \sum_{\mu \in Spec_H(\pi_\lambda^G)} m^H(\lambda ,\mu) \,\, \pi_\mu^H $$ Equality in the previous formula means unitary equivalence. On the right hand side, sum,  means  the Hilbert sum of the family of representations of $H$ pointed out. For a proof c.f. \cite{kobtwo}. \\
 (0.2a) In  \cite{gw}, \cite{oshi} it is shown that every element of $Spec_H(\pi_\lambda^G)$ is dominant with respect to $\Psi_{H,\lambda}.$   Moreover, since we are dealing with square integrable representations, we have that $\pi_\lambda^G$ has an admissible restriction to $H$ if and only if $\pi_\lambda^G$ is an admissible representation of $H_0.$ This follows from \cite{dv}. \\

 To each system of positive roots  $\Psi_\lambda,$ in \cite{dv}, we have attached a connected normal subgroup    $K_1 := Z_1(\Psi_\lambda)K_1(\Psi_\lambda)$ of $K$ and we have shown:

\smallskip
\noindent
 (0.3) $\pi_\lambda^G$ restricted to $H$ is admissible if and only if $Z_1(\Psi_\lambda)K_1(\Psi_\lambda)$  is a subgroup of $L.$

  In section 1 we list the 5-tuples  $(G, H, H_0, \Psi_\lambda, Z_1(\Psi_\lambda)K_1(\Psi_\lambda)) $ so that $G$ is simple and  $\pi_\lambda^G$ has an admissible restriction to $H.$ In the course of the note we list the systems $\Psi_{H,\lambda}, \Psi_{H_0,\lambda}.$

\smallskip
\noindent
(0.4) Let $\k_2 =\k_2(\Psi_\lambda)$ denote the complementary ideal to $\k_1.$  Thus, we have a decomposition, as a direct  sum of  ideals, $$\k =\k_1 + \k_2 , \,\, \t =\t_1 + \t_2, \,\,\,  \t^\star = \t_1^\star  + \t_2^\star. $$  as we  well as the product of normal subgroups $K=K_1 K_2.$ Owing to (0.3),  $K_1$ is a  subgroup of $L$ hence  we have the decomposition $L=K_1 (L\cap K_2)$ and the decomposition $\u=\t_1 + (\l \cap \t_2).$ For a reductive Lie algebra $\s$ we write  $\s=\z_\s +\s_{ss}$ where $\z_\s$ denotes the center and $\s_{ss}$ denotes the semisimple factor. For some examples, the decomposition $K=K_1 K_2, \, T=T_1 T_2$ and the corresponding decomposition for $\u$ are not completely compatible with the decompositions $\h=\z_\h + \h_{ss}, \l=\z_\l +\l_{ss} = \z_h  + (\z_\l \cap \h_{ss}) + \l_{ss}, \z_\l =\z_\h + \z_\l \cap \h_{ss}.$  However, this is not a drawback for the statements to follow, we dealt with this matter in (1.6).   For $\lambda \in \mathfrak t^\star, $ we write $\lambda =\lambda_1 + \lambda_2 =(\lambda_1, \lambda_2),$  with $\lambda_i \in \t_i^\star, i=1,2.$   For $\gamma \in \u^\star, $ we write $\gamma = \gamma_1 + \gamma_2=(\gamma_1, \gamma_2)$ with $\gamma_1 \in \t_1^\star, \gamma_2 \in (\l \cap \t_2)^\star.$ Because of our hypothesis on $\lambda,$ it follows that $\lambda_2$ is Harish-Chandra parameter for  $K_2,$ (c.f. 1.46), hence, there exists an representation $  \pi_{\lambda_2}^{K_2} $ of infinitesimal character $\lambda_2.$   The branching law for the restriction of  $  \pi_{\lambda_2}^{K_2} $ to $L\cap K_2$ is $$ res_{L\cap K_2}( \pi_{\lambda_2}^{K_2}) = \sum_{\nu_2 \in Spec_{L\cap K_2}(
 \pi_{\lambda_2}^{K_2})} m^{K_2, L\cap K_2}(\lambda_2 , \nu_2) \,\, \pi_{\nu_2}^{L\cap K_2}.$$

For an irreducible square integrable representation $\pi_{\mu}^{H_0}$ of $H_0.$
      The assumption, $L$ is a maximal compact subgroup of $H_0, $ gives rise to the decomposition $$res_L (\pi_{\mu}^{H_0})= \sum_{\nu \in Spec_L(\pi_{\mu}^{H_0}) } m^{H_0,L}(\mu, \nu)\,\, \pi_\nu^L $$ Here, as before, $\pi_\nu^L $ denotes de irreducible representation of $L$ of infinitesimal character $\nu$ dominant with respect to $\Delta_0,$ and
 $Spec_L(\pi_{\mu}^{H_0})$  denotes the set of Harish-Chandra parameters  of irreducible representations of $L$ that occurs in $res_L (\pi_{\mu}^{H_0}).$ In section 1 we study when   for  $\nu_2 \in Spec_{L\cap K_2}(\pi_{\lambda_2}^{K_2})  $ the pair  $(\lambda_1, \nu_2)$ is a Harish-Chandra parameter for $H_0.$  We first  show $(\lambda_1, \nu_2)$ is regular and dominant for $\Psi_{H_0, \lambda} .$  However, $e^{(\lambda_1, \nu_2)+\rho_{H_0}}$ may not be a character of $U.$ This is not a problem because  twice a Harish-Chandra parameter always   lifts to a character. Thus,  $e^{(\lambda_1, \nu_2)+\rho_{H_0}}$ is a character of a two-fold cover of $U.$ Hence, as usual, we replace $G$ by a two fold cover and then we  have that $(\lambda_1, \nu_2)$ is a Harish-Chandra parameter. This assumption has as a consequence, that any Harish-Chandra parameters for $H$ which occurs in $res_H(\pi_\lambda^G) $ is also a Harish-Chandra parameters for $L.$ Now, we may state the first result of this note.

\begin{thm}Assume $\pi_\lambda^G$ has an admissible restriction to $H.$ Then, we have
 $$ m^{H}(\lambda, \mu)=\sum_{ \nu_2 \in Spec_{L\cap K_2}(
\pi_{\xi_2}^{K_2})} m^{H_0, L}((\lambda_1, \nu_2), \mu)  \,\, m^{K_2, L\cap K_2}(\lambda_2, \nu_2)$$
for each $\mu \in Spec_H(\pi_\lambda^G).$
\end{thm}

  The lowest $K-$type $\pi_\xi^K$  of $\pi_\lambda^G$, \cite{vothesis}, decomposes as the outer tensor product of irreducible representations $$\pi_\xi^K = \pi_{\xi_1}^{K_1} \boxtimes \pi_{\xi_2}^{K_2}$$

  It is known, \cite{gw}, it may happens that $\pi_{\xi_2}^{K_2}$ is the trivial representation of $K_2.$ Under this hypothesis, Bent Orsted and Birgit Speh in \cite{os}, \cite{Speh}  conjectured that for a convenient Harish-Chandra parameter $\lambda^\prime$  for $H_0,$   the set  of Harish-Chandra parameters of the  $L-$types of $res_L (\pi_{\lambda^\prime}^{H_0})$ is equal to the set of Harish-Chandra parameters of $res_H(\pi_\lambda^G)$ and that  the multiplicity functions $m^H(\lambda, ?) , m^{H_0, L}(\lambda^\prime, ?)$ are equal. Actually, their conjecture is stated for a family of Zuckerman modules $A_\mathfrak q (\lambda)$ which includes the family of square integrable  representations so that $\pi_{\xi_2}^{K_2}$ is the trivial representation. As a consequence of Theorem 1 we have that their conjecture is true. In fact,

\begin{thm}  Assume $\pi_\lambda^G$ is an admissible representation of $H$ as well as that the lowest $K-$type of $\pi_\lambda^G$ is an irreducible representation of $K_1.$  Then,  $$ Spec_H(\pi_\lambda^G)= Spec_L(\pi_{(\lambda_1, \rho_{L\cap K_2})}^{H_0}), \,\,\,\, m^H(\lambda, \mu)=m^{H_0, L}((\lambda_1, \rho_{L\cap K_2}), \mu)$$ for every $\mu \in Spec_L(\pi_{\lambda }^{H_0}).$
\end{thm}

One consequence of Theorem 2 and Lemma 2.12 in \cite{schannals} is.
\begin{cor} Every $\mu \in Spec_H(\pi_\lambda^G)$  is equal to $(\lambda_1 , \rho_{L\cap K_2}) + \rho_0 +B $ where $B$ is a sum of roots in $\Psi_{H_0, \lambda} \cap \Phi(\p \cap \mathfrak q, \u)= \Psi_{H_0,(\lambda_1, \rho_{L\cap K_2})} \cap \Phi_n (\h_0, \u)$ and $\rho_0$ is equal to one half of the sum of the roots in $\Psi_{H_0,(\lambda_1, \rho_{L\cap K_2})} \cap \Phi_n (\h_0, \u).$
\end{cor}

Theorem 1 and Theorem 2 are somewhat  analogous to results of Gross-Wallach \cite{gw} and statements on branching laws presented in the book of Knapp \cite{knapp}. Besides,  Theorem 2 has a resemblance to results of Kobayashi proved in \cite{kobmem}.\\ Theorem 1 and Theorem 2 coupled with the work of Baldoni Silva-Vergne, \cite{bsv},  provide an effective method to compute the  multiplicity   and the Harish-Chandra parameter for the irreducible factors in  the decomposition  of $res_H (\pi_\lambda^G)$  as $H-$module.

In the fourth section of this note we analyze the subspace $\mathcal L_\lambda$  spanned by the lowest $L-$type of the totality irreducible $H-$factors and for a scalar holomorphic discrete series we obtain a more precise description of  $\mathcal L_\lambda.$  We also consider   the relation between $Hom_H(\pi_\mu^H, \pi_\lambda^G)$ and the space of intertwining operators from the lowest $L-$type of $\pi_\mu^H$ into $res_L(\pi_\lambda^G).$ In the fifth  section we show, for outer tensor product of two holomorphic discrete series representations, a similar result to Theorem 1 as well as a similar result to the one obtained in section 4. In section 6, we collect notation and complete some case by case proofs.

\section{Description of the subgroup $K_1.$}

We maintain the notation of the previous section.  For each root $\alpha \in \Phi(\g, \t),$ let $\g_\alpha$ denote the root space associated to $\alpha,$  then, by definition, $\k_1(\Psi_\lambda)$ is the ideal of $\k$  spanned by $[\mathfrak g_\alpha, \mathfrak g_\beta], \alpha, \beta \in  \Psi_\lambda \cap  \Psi_n. $ We define $K_1(\Psi_\lambda)$ the analytic subgroup of $K$ which corresponds to $\k_1(\Psi_\lambda).$   Thus, $K_1(\Psi_\lambda)=\{ 1\}$ if and only if $\Psi_\lambda$ is a holomorphic system. With respect to $Z_1(\Psi_\lambda),$    whenever $G/K$ is not an hermitian symmetric space, $Z_1(\Psi_\lambda)$ is defined to be equal to the trivial group. For a holomorphic system $\Psi_\lambda$  we define $Z_1(\Psi_\lambda)$ to be equal to the identity connected component of the center of $K,$ hence,  $Z_1(\Psi_\lambda)$ is a one dimensional torus.  We claim:\\

 (1.1) Assume $G/K$ is an Hermitian symmetric space and let $\Psi_\lambda$   be a nonholomorphic system of positive roots such that $\pi_\lambda^G$ is an admissible representation of $H,$ then the subgroup  $Z_1(\Psi_\lambda)$ defined in \cite{dv} is  the trivial group.

\smallskip
\noindent
  To verify the claim we notice that  the hypothesis on admissibility together with the tables in \cite{kosh}, or else the computations in \cite{dv}, \cite{dvc}, yields,  \\
  (1.2)  If   $G/K$ is an Hermitian symmetric space and $\pi_\lambda^G$ is not a holomorphic representation with admissible restriction to $H,$   then  the pair $(\g_\mathbb R, \h_\mathbb R)$  is one of \begin{center} $ (\mathfrak {su}(m,n), \s (\mathfrak u(m,k) +\u (n-k))), \, (\mathfrak {su}(2,2n), \mathfrak {sp}(1,n)),$  \\ $\,(\mathfrak {so}(2m,2), \mathfrak {so}(2m,1)), \,   \, (\mathfrak {su}(2,2), \mathfrak {sp}(1,1)).$ \end{center} In the next paragraph, for each of these pairs  $(G,H),$  we write compact Cartan subgroups and list the  systems of positive roots $\Psi_\lambda$ so that $res_H (\pi_\lambda^G)$ is an admissible representation,   for  each case, we verify $Z_1(\Psi_\lambda)=\{0\}$ and we compute $\Psi_{H_0, \lambda}, \Psi_{H, \lambda}.$

   \noindent
    $\mathbf{I-1, I-2}.$   $(\mathfrak {su}(m,n), \s (\mathfrak u(m,k) +\u (n-k))).$ We fix as Cartan subalgebra   $\t$ of $\mathfrak {su}(m,n)$ the set of diagonal matrices in $\mathfrak {su}(m,n).$   For certain orthogonal basis $\epsilon_1, \dots ,\epsilon_p, \delta_1, \dots, \delta_q$ of the dual vector space to the subspace of diagonal matrices in $\mathfrak{gl}(m+n, \mathbb C),$ we may, and will choose  $\Delta =\{ \epsilon_r - \epsilon_s, \delta_p -\delta_q, 1 \leq  r < s \leq m, 1 \leq p < q \leq n \},$ the set of noncompact roots is $ \Phi_n= \{ \pm (\epsilon_r - \delta_q) \}.$ We recall the  positive roots systems for $\Phi(\g, \t)$  containing $\Delta$ are in a bijective  correspondence with the totality of lexicographic orders for the basis $\epsilon_1, \dots ,\epsilon_m, \delta_1, \dots, \delta_n$ which contains the "suborder" $\epsilon_1 > \dots > \epsilon_m, \,\, \delta_1 > \dots > \delta_n.$ The two holomorphic systems correspond to the orders $\epsilon_1 > \dots > \epsilon_m > \delta_1 > \dots > \delta_n ; \,\,  \delta_1 > \dots > \delta_n >\epsilon_1 > \dots > \epsilon_m.$ We fix $1 \leq a \leq m-1, $  in \cite{dvc} is verified that for the system of positive roots $\Psi_a$ corresponding to the order $\epsilon_1> \dots > \epsilon_a > \delta_1> \dots> \delta_n > \epsilon_{a+1}  > \dots > \epsilon_m , $  we have $\k_1(\Psi_a)= \mathfrak{su}(m).$ We fix $1 \leq b \leq n-1 $ and let $\tilde{\Psi}_b$ denote the set of positive roots associated to the order $\delta_1 >  \dots > \delta_b >\epsilon_1 > \dots > \epsilon_m > \delta_{b+1} > \dots > \delta_n,$ then $\k_1(\tilde{\Psi}_b)=\mathfrak{su}(n).$ For any other nonholomorphic system $\Psi$  we have $\k_1(\Psi)= \mathfrak{su}(m)+ \mathfrak{su}(n).$ Thus, (0.3) forces $\Psi_\lambda$ to be equal to either $\Psi_a$ or $\tilde{\Psi}_b.$   A direct computation \cite{dvc} verifies $\mathbb R^+ (\Psi_a \cap \Phi_n) \cap i\mathfrak (\z_\k)_\mathbb R^\star = \mathbb R^+ (\tilde{\Psi }_b \cap \Phi_n) \cap i\mathfrak (\z_\k)_\mathbb R^\star =\{0\}, $ hence the definition of $Z_1( \Psi_\lambda) $ implies $Z_1( \Psi_\lambda) =\{1 \}$. Thus, we have verified  claim (1.1) for the pair $(\mathfrak {su}(m,n), \mathfrak s(\u(m,k)+ \u (n-k))).$    The root systems for  $(\mathfrak h, \mathfrak t)$ and its dual are:
       \begin{multline*}
     \Phi(\mathfrak h, \mathfrak t)=\{ \pm (\epsilon_r -\epsilon_s), \pm(\delta_p -\delta_q), \pm (\epsilon_i -\delta_j),
    1 \leq r < s \leq m,  \\  1 \leq p <  q \leq k, \,\,   or, \,\, k+1 \leq p <  q \leq n, 1\leq i \leq m,  1 \leq j \leq k \}.
    \end{multline*}
   \begin{multline*} \Phi(\mathfrak h_0, \mathfrak t)=\{ \pm (\epsilon_r -\epsilon_s), \pm(\delta_p -\delta_q), \pm (\epsilon_i -\delta_j),  1 \leq r < s \leq m, \\  1 \leq p <  q \leq k\,\, or\,\, k+1 \leq p < q \leq n,  1 \leq i \leq m,  k+1 \leq j \leq n \}.
   \end{multline*}

       The system $\Psi_{H,\lambda},$  $\Psi_{H_0,\lambda}$ which correspond to  $\Psi_a$ are the system associated to the respective lexicographic orders  $$ \epsilon_1> \dots > \epsilon_a > \delta_1> \dots> \delta_k > \epsilon_{a+1}  > \dots > \epsilon_m $$ $$\epsilon_1> \dots > \epsilon_a > \delta_{k+1}> \dots> \delta_n > \epsilon_{a+1}  > \dots > \epsilon_m .$$ For $\tilde{\Psi_b}$ the description of $\Psi_{H,\lambda},$  $\Psi_{H_0,\lambda}$ is similar.

         From now on, $q_\u $ denotes the restriction map from $\t^\star$ onto $\u^\star.$

   \smallskip
\noindent
   $\mathbf{II-1}. $
     $(\mathfrak {su}(2,2n), \mathfrak {sp}(1,n)), n \geq 1 .$ We use the notation in $\mathbf{A1}$. The automorphism $\sigma$ acts on $\t^\star$ as    $\sigma (\epsilon_1) =-\epsilon_2,   \sigma (\delta_i )=  -\delta_{2n-i +1}, i=1, \dots, 2n .$  Hence, a basis of  $\u^\star $ is  $q_\mathfrak u (\epsilon_1 -\epsilon_2), q_\mathfrak u (\delta_r -\delta_{2n-r+1}), r=1, \dots n . $
   Because of the way $\sigma$ acts on $\t^\star$ we have that $\beta \not= \sigma (\beta)$ for any  noncompact root. For each root $\beta,$ we fix $0 \not= Y_\beta \in \mathfrak g_\beta.$ Then, $\p \cap \mathfrak q = \sum_{\beta \in \Phi_n} \mathbb C (Y_\beta - \sigma Y_\beta ) $ and $ \h \cap \p = \sum_{\beta \in \Phi_n} \mathbb C (Y_\beta + \sigma Y
   _\beta ).$ Hence, $\Phi_n (\h , \u)=\Phi_n(\h_0, \u)=\{\pm q_\u (\epsilon_1 -\delta_i), i=1, \dots 2n \}.$ Certainly, $\Phi_c (\h, \u)=\Phi_c (\h_0, \u)= \{ \pm(\epsilon_1 -\epsilon_2) , q_\mathfrak u (\delta_i -\delta_j), i\not= j \}.$  The unique  possibility for $\k_1(\Psi_\lambda)$ to be contained in $\mathfrak{sp}(1,n)$ is for  $\Psi_\lambda = \Psi_1.$    The simple roots for $\Psi_1$ are $\epsilon_1 -\delta_1, \delta_1 -\delta_2, \dots , \delta_{2n-1} -\delta_{2n}, \delta_{2n} -\epsilon_2. $ We notice  $\sigma (\Psi_1)=\Psi_1$ and $\k_1(\Psi_\lambda) =\mathfrak{su}
   _2( \epsilon_1 -\epsilon_2).$   We define  $\Psi_{H, \lambda} = \Psi_{H_0 , \lambda} $
  and we set   $\Psi_{H,\lambda} \subset \Phi(\h, \u)$ to be the system of positive roots for the simple roots $$\{\alpha_1 := q_\u(\epsilon_1 -\delta_1), \alpha_{i+1}:= q_\u( \delta_i -\delta_{i+1}), i=1, \dots, n-1, \alpha_{n+1}:= (\delta_n -\delta_{n+1}) \}. $$  For the pair $(SU(2,2), Sp(1,1)),$ we  have the extra possibility  $\Psi_\lambda= \tilde{\Psi}_1.$ In \cite{dvc} we verify  for any of these systems   $Z_1(\Psi_\lambda)$ is equal to the trivial group.

   \smallskip
\noindent
    $\mathbf{II-5} $  $ (\mathfrak {so}(2m,2), \mathfrak {so}( 2m,1)), m \geq 2.$ We  choose an orthogonal basis $\{ \epsilon_1,   \dots , \epsilon_m, \delta_1 \}$ of $i\t_\mathbb R^\star$    so that $\Delta =\{ \epsilon_k \pm \epsilon_s, 1 \leq k < s \leq m \}, \Phi_n =\{ \pm (\epsilon_j \pm \delta_1), 1 \leq j \leq m \}.$ In this case,  $\mathfrak z_K^\star =\mathbb C \delta_1. $ The systems of positive roots $\Psi_\lambda$ containing $\Delta$ are parameterized by the lexicographic orders  $S_{\pm a} := \epsilon_1 > \dots > \epsilon _a >\pm \delta_1 > \epsilon_{a+1} > \dots > \epsilon_{m-1} > \pm \epsilon_m $ for $a=0, 1, \dots, m.$ The two holomorphic systems corresponds to the parameters $\pm 0.$   In this case $\sigma$ acts on $\t^\star$ by  the rule  $\sigma (\delta_1)=-\delta_1, \sigma (\epsilon_j)=\epsilon_j, j=1, \dots, m.$ Hence,  $\u^\star$ is spanned by $\epsilon_1, \dots, \epsilon_m. $ It readily follows that $\sigma$ leaves invariant the system of positive roots associated to $ S_{\pm m}.$  Also, $\sigma (\beta) \not= \beta$ for every noncompact root. Hence, as in the previous case, $\Phi_n(\h, \u)=\Phi_n (\h_0, \u).$ Here, $$ K_1(\Psi_\lambda)= \left\{ \begin{array}{ll} \mathfrak{so}(2m) & \mbox{if} \,\, m>2, \,\,\mbox{  $\Psi_\lambda$ nonholomorphic } \\ \mathfrak{su}_2(\epsilon_1 \pm \epsilon_2) & \mbox{if} \,\, m=2,\,\,\,\, \Psi_\lambda \leftrightarrow S_{\pm 2} \\ \mathfrak{so}(4) & \mbox{if} \,\, m=2, \,\,\,\, \Psi_\lambda \leftrightarrow S_{\pm 1} \end{array} \right. $$
    From  actual computation or the computations in \cite{dvc} we have $\mathbb R^+ \Psi_\lambda \cap \mathfrak z_K^\star =\{0\}$ if and only if $\Psi_\lambda$ is  associated to  one of the orders $ S_{\pm m}.$   From the tables in \cite{kosh}, or  computations in \cite{dv}, it follows that $\pi_\lambda^G$ has an admissible restriction to $H$ if and only if  $\Psi_\lambda$ is  associated to  one of the orders $ S_{\pm m}.$   Hence, $Z_1(\Psi_\lambda)=\{1 \}$ for both cases.     Let $\tau $ denote $\pm.$ For the system of positive roots $\Psi_\tau $ associated to the order $S_{\tau m}$ we define   $\Psi_{H, \lambda}=\Psi_{H_0, \lambda}$ to be the system of positive roots for the simple roots $$ (\epsilon_j -\epsilon_{j+1}), (\epsilon_{m-1} -\tau  \epsilon_m), q_\u (\tau \epsilon_m -\delta_1)=\tau \epsilon_m, \,\, j=1, \dots, m-2.$$  This concludes the verification of (1,1). \begin{flushright} $
 \Box  $ \end{flushright}

 For an arbitrary symmetric pair $(G,H),$ whenever $\pi_\lambda^G$ is an admissible representation of $H,$  because of (1.1), we define,
  $$ K_1= \left\{ \begin{array}{ll} Z_K & \mbox{\, if $\Psi_\lambda$ holomorphic } \\ K_1(\Psi_\lambda) & \phantom{x} \mbox{otherwise} \end{array} \right. $$
In the next tables we present the 5-tuple so that: $(G, H)$ is a symmetric pair, $H_0$ is the associated group to $H,$ $\Psi_\lambda$ is a  system of positive such that $\pi_\lambda^G$ is an admissible representation of $H,$ and $K_1 =Z_1( \Psi_\lambda) K_1(\Psi_\lambda).$  Actually, instead of writing  Lie groups we write their respective Lie algebras. Each  table is in part a reproduction of tables in \cite{kosh}. The tables  can also be computed by means of the techniques presented in \cite{dv}. Note that each table is "symmetric" when we replace $H$ by $H_0.$ As usual, $\alpha_m$ denotes the highest root in $\Psi_\lambda$ and $\mathfrak{su}_2(\alpha)$ denotes the compact real  subalgebra of $\mathfrak g_\mathbb R$  spanned by the root vectors corresponding to the compact root $\alpha.$  In this and other sections we complete the notation for the objects in the tables.

{\tiny

\begin{center}
\begin{tabular}{|c| c | c| c |c |} \hline  $G$  &  $H$  & $H_0$  &  $\Psi_\lambda$  &  $ K_1$  \\
\hline  $\mathfrak{su}(m,n)$  &  $\mathfrak{su}(m,k)\oplus \mathfrak{su}(n-k)\oplus \mathfrak{u}(1)$  &  $\mathfrak{su}(m,n-k)\oplus \mathfrak{su}(k)\oplus \mathfrak{u}(1)$   &  $\Psi_a$  &  $ \mathfrak{su}(m)$  \\
\hline  $\mathfrak{su}(m,n)$  &  $\mathfrak{su}(k,n)\oplus \mathfrak{su}(m-k)\oplus \mathfrak{u}(1)$  &  $\mathfrak{su}(m-k,n)\oplus \mathfrak{su}(k)\oplus \mathfrak{u}(1)$   &  $\tilde{\Psi_b}$  &  $ \mathfrak{su}(n)$  \\
\hline   $\mathfrak{so}(2m,2n)$   &   $\mathfrak{so}(2m,2k)\oplus \mathfrak{so}(2n-2k)$ & $\mathfrak{so}(2m,2n-2k)\oplus \mathfrak{so}(2k)$   &   $\Psi_{\pm}$  & $ \mathfrak{so}(2m)$   \\
\hline   $\mathfrak{so}(2m,2n+1)$   &   $\mathfrak{so}(2m,k)\oplus \mathfrak{so}(2n+1-k)$ &   $\mathfrak{so}(2m,2n+1-k)\oplus \mathfrak{so}(k)$   &   $\Psi_{\pm}$   &   $ \mathfrak{so}(2m)$   \\
\hline   $\mathfrak{so}(4,2n),  n>2$   &   $\mathfrak{u}(2,n)_1$   &   $w\mathfrak{u}(2,n)_1$   &   $\Psi_{1\, -1}$   &   $ \mathfrak{su}_2(\alpha_m)$   \\
\hline   $\mathfrak{so}(4,2n),  n>2$   &   $\mathfrak{u}(2,n)_2$   &   $w\mathfrak{u}(2,n)_2$   &   $\Psi_{1 \, 1}$   &   $ \mathfrak{su}_2(\alpha_m)$   \\
\hline   $\mathfrak{so}(4,4)$   &   $\mathfrak{u}(2,2)_{1\,1}$   &   $w\mathfrak{u}(2,2)_{11}$   &   $\Psi_{1\, -1},  \, w_{\epsilon,\delta}\Psi_{1 \, -1}$   &   $ \mathfrak{su}_2(\alpha_m)$ \\
\hline
$\mathfrak{so}(4,4)$   &   $\mathfrak{u}(2,2)_{12}$   &   $w\mathfrak{u}(2,2)_{12}$   &   $\Psi_{1\, -1},  \, w_{\epsilon,\delta}\Psi_{1\, 1}$   &    $ \mathfrak{su}_2(\alpha_m)$ \\
\hline   $\mathfrak{so}(4,4)$   &   $\mathfrak{u}(2,2)_{21}$   &   $w\mathfrak{u}(2,2)_{21}$   &   $\Psi_{1\, 1}, \,  w_{\epsilon,\delta}\Psi_{1\, -1}$   &   $ \mathfrak{su}_2(\alpha_m)$ \\
\hline   $\mathfrak{so}(4,4)$   &   $\mathfrak{u}(2,2)_{22}$   &   $w\mathfrak{u}(2,2)_{22}$   &   $\Psi_{1\,1}, \, w_{\epsilon,\delta}\Psi_{1\,1}$   &   $ \mathfrak{su}_2(\alpha_m)$ \\
\hline
  $\mathfrak{sp}(m,n)$   &   $\mathfrak{sp}(m,k)\oplus \mathfrak{sp}(n-k)$   &   $\mathfrak{sp}(m,n-k)\oplus \mathfrak{sp}(k)$   &   $\Psi_+$   &   $ \mathfrak{sp}(m)$   \\
\hline   $\mathfrak f_{4(4)}$   &   $\mathfrak{sp}(1,2)\oplus \mathfrak{su}(2)$ & $\mathfrak{so}(5,4)$   &   $\Psi_{BS}$   &   $ \mathfrak{su}_2(\alpha_m)$   \\
\hline
\phantom{ } $\mathfrak{e}_{6(2)}$   &   $\mathfrak{so}(6,4)\oplus \mathfrak{so}(2)$ & $\mathfrak{su}(4,2)\oplus \mathfrak{su}(2)$   &   $\Psi_{BS}$   &   $ \mathfrak{su}_2(\alpha_m)$   \\
\hline
  $\mathfrak{e}_{7(-5)}$   &   $ \mathfrak{so}(8,4)\oplus \mathfrak{su}(2)$ &   $\mathfrak{so}(8,4)\oplus \mathfrak{su}(2)$   &   $\Psi_{BS}$ & $ \mathfrak{su}_2(\alpha_m)$   \\
\hline
  $\mathfrak{e}_{7(-5)}$   &   $\mathfrak{su}(6,2)$   &   $\mathfrak{e}_{6(2)}\oplus \mathfrak{so}(2)$   &   $\Psi_{BS}$   &   $ \mathfrak{su}_2(\alpha_m)$   \\
\hline
  $\mathfrak{e}_{8(-24)}$   &   $\mathfrak{so}(12,4)$   &   $\mathfrak{e}_{7(-5)}\oplus \mathfrak{su}(2)$   &   $\Psi_{BS}$   &   $ \mathfrak{su}_2(\alpha_m)$   \\
\hline
\end{tabular} \\
Table I. Case $U=T, \Psi_\lambda$ nonholomorphic

\end{center}

\begin{center}
\begin{tabular}{|c| c | c| c |c |} \hline  $G$    &  $H$  &  $H_0$  & $\Psi_\lambda$  &  $ K_1$  \\
\hline  $\mathfrak {su}(2,2n), \, n> 2$  & $\mathfrak {sp}(1,n)$  &  $\mathfrak {sp}(1,n)$  &  $\Psi_1$  &  $ \mathfrak{su}_2(\alpha_m)$  \\
\hline  $\mathfrak{su}(2,2)$  & \phantom{xx} $\mathfrak{sp}(1,1)$ &  $\mathfrak{sp}(1,1)$ &  $\Psi_1$  &  $ \mathfrak{su}_2(\alpha_m)$  \\
\hline
   $\mathfrak{su}(2,2)$ &  $\mathfrak{sp}(1,1)$  &  $\mathfrak{sp}(1,1)$  &  $\tilde{\Psi}_1$  & $ \mathfrak{su}_2(\alpha_m)$  \\
\hline
  $\mathfrak{so}(2m,2n), m>1$  &  $ \mathfrak{so}(2m,2k+1) + \mathfrak{so}(2n-2k-1)$  &  $\mathfrak{so}(2m,2n-2k-1)+ \mathfrak{so}(2k+1)$  &  $\Psi_{\pm}$  &  $ \mathfrak{so}(2m)$ \\
  \hline
  $\mathfrak{so}(2m,2), m>2 $  &  $ \mathfrak{so}(2m,1) $  &  $\mathfrak{so}(2m,1)$  &  $\Psi_{\pm}$  &  $ \mathfrak{so}(2m)$ \\
\hline
 $\mathfrak{e}_{6(2)}$  &  $\mathfrak{f}_{4(4)}$  &  $\mathfrak{sp}(3,1)$  &  $\Psi_{BS}$  &  $ \mathfrak{su}_2(\alpha_m)$  \\
\hline
\end{tabular} \\

Table II, Case $U\not=T, \Psi_\lambda$ non holomorphic
\end{center}

\begin{center}
\begin{tabular}{|c| c | c|}
\hline   $G$  &   $H \, $ (a)  &  $H_0 \,   $  (b)   \\
 \hline   $\mathfrak{su}(m,n), m\not= n$   &     $\mathfrak{su}(k,l)+\mathfrak{su}(m-k,n-l)+ \u(1)$ & $\mathfrak{su}(k,n-l)+\mathfrak{su}(m-k,l)+ \u(1)$    \\
 \hline   $\mathfrak{su}(n,n)$  & $\mathfrak{su}(k,l)+ \mathfrak{su}(n-k,n-l)+\u(1)$ &$\mathfrak{su}(k,n-l)+ \mathfrak{su}(n-k,l)+\u(1)$       \\
 \hline   $\mathfrak{so}(2,2n)$  & $\mathfrak{so}(2,2k)+ \mathfrak{so}(2n-2k)$   & $\mathfrak{so}(2,2n-2k)+ \mathfrak{so}(2k)$      \\
\hline   $\mathfrak{so}(2,2
n)$  & $\mathfrak{u}(1,n)$ & $ \mathfrak{u}(1,n)$        \\
\hline   $\mathfrak{so}(2,2n+1)$  & $\mathfrak{so}(2,k)+ \mathfrak{so}(2n+1-k)$ & $\mathfrak{so}(2,2n+1-k)+ \mathfrak{so}(k)$       \\
\hline   $\mathfrak{so}^\star (2n)$  & $\mathfrak{u}(m,n-m)$ & $\mathfrak{so}^\star(2m)+ \mathfrak{so}^\star(2n-2m)$       \\
\hline   $\mathfrak{sp}(n, \mathbb R)$  & $\mathfrak{u}(m,n-m)$ & $\mathfrak{sp}(m, \mathbb R)+ \mathfrak{sp}(n-m, \mathbb R)$       \\
\hline $\e_{6(-14)}$ & $\mathfrak{so}(2,8)+\mathfrak{so}(2) $ & $\mathfrak{so}(2,8)+\mathfrak{so}(2)$ \\
\hline $\e_{6(-14)}$ & $\mathfrak{su}(2,4)+\mathfrak{su}(2) $ & $\mathfrak{su}(2,4)+\mathfrak{su}(2)$ \\
\hline $\e_{6(-14)}$ & $\mathfrak{so}^\star(10)+\mathfrak{so}(2) $ & $\mathfrak{su}(5,1)+\mathfrak{sl}(2, \mathbb R)$ \\
\hline $\e_{7(-25)}$ & $\mathfrak{so}^\star(12)+\mathfrak{su}(2) $ & $\mathfrak{su}(6,2)$ \\
\hline $\e_{7(-25)}$ & $\mathfrak{so}(2,10)+\mathfrak{sl}(2, \mathbb R) $ & $\e_{6(-14)} + \mathfrak{so}(2)$ \\
\hline   $\mathfrak{su}(n,n)$  & $\mathfrak{so}^\star(2n)$ &$\mathfrak{sp}(n,\mathbb R)$       \\
\hline   $\mathfrak{so}(2,2n)$  & $\mathfrak{so}(2,2k+1)+ \mathfrak{so}(2n-2k-1)$   & $\mathfrak{so}(2,2n-2k-1)+ \mathfrak{so}(2k+1)$      \\
\hline
\end{tabular} \\
Table III, $\pi_\lambda^G $ holomorphic Discrete Series. \\
The last two lines show the unique holomorphic pairs so that $U \not= T.$
\end{center}
}

For an arbitrary parameter of Harish-Chandra, $\lambda,$ we recall  that we denote the complementary ideal to $\k_1$  by $\k_2.$ Hence, the root system $\Phi (\k_2, \t_2)$ is a subsystem of $\Phi_c. $ In \cite{dv} we find a proof of: \\ (1.3) Any simple root for $\Delta \cap \Phi (\k_2, \t_2)$ is a simple root for $\Psi_\lambda.$

   Let $W_S$ denote the Weyl group of a compact connected Lie group $S.$  Then, (1.3) and a result in  \cite{hs} yields: \\ (1.4) for $w \in W_{K_2},$ we have $w (\Psi_\lambda \cap \Phi_n) = \Psi_\lambda \cap \Phi_n.$ \\
 We define, $$
\rho_c :=\frac{1}{2} \sum_{\alpha \in \Delta}\alpha, \,\,  \rho_n^\lambda := \frac{1}{2} \sum_{\beta \in \Psi_\lambda \cap \Phi_n }\beta,  \,\, \rho := \rho_c + \rho_n^\lambda, \,\, \rho_{K_2}  :=\frac{1}{2} \sum_{\alpha \in \Delta \cap \Phi(\k_2, \t_2)}\alpha $$
  (1.4a), $\rho_n^\lambda \in \t_1^\star, $  because for  $\alpha$ a compact simple root for $\Psi_\lambda,$ we have $(\alpha, \rho_n^\lambda)=0.$ As before, we write $\lambda =(\lambda_1, \lambda_2),$ we have: \\ (1.4b) $\lambda_2$ is a Harish-Chandra parameter for $K_2.$  Indeed,   the character $e^{\lambda + \rho} $ restricted to $T \cap K_2$ is equal to the function $e^{\lambda_2 + \rho_{K_2}}.$ \\ (1.5) We now verify that for $\nu_2 \in Spec_{L\cap K_2}(\pi_{\lambda_2}^{K_2}) ,$ then the pair $ (\lambda_1, \nu_2)$ is either a Harish-Chandra parameter for $H_0$ or is a Harish-Chandra parameter for a two fold cover of $H_0.$ Always, $(\lambda_1, \nu_2)$ is dominant and regular for $\Psi_{H_0, \lambda}.$

\noindent
    Indeed,  sometimes, $\rho_G -\rho_{H_0}$ may not lift to a character of $U,$ hence $(\lambda_1 , \nu_2) + \rho_{H_0}$ may not lift to a character of $U,$ however, twice of $ (\lambda_1 , \nu_2) + \rho_{H_0}$ does lift to a character of $U.$ Thus, we only need to verify $ (\lambda, \nu_2) >0 $ for  any noncompact root in $\Psi_{H_0, \lambda}.$ We first analyze the case $U=T.$ Owing to a theorem of Kostant, any Harish-Chandra parameter of $res_{L\cap K_2} (\pi_{\lambda_2}^{K_2}) $ belongs to the convex  hull of the set $s \lambda_2 , s \in W_{K_2}. $ Whence, $\nu_2 =\sum_{s \in W_{K_2}}c_{s} s\lambda_2,$ with $c_{s} \geq 0, \sum c_{s}=1.$  Because of (1.4), for $\beta \in \Psi_\lambda \cap \Phi_n$  we have $((\lambda_1, \nu_2), \beta)=\sum_{s \in W_{K_2}} c_{s}((\lambda_1, \lambda_2), s^{-1}\beta) >0.$  Thus, $(\lambda_1, \nu_2)$ is a Harish-Chandra parameter for either $H_0$ or a two fold cover of $H_0.$   The case $U \not=T $ follows from the previous computation and the observation that $H_0$ has only one  noncompact simple factor, hence, the inner product defined  by the respective Killing form's in the Cartan subalgebras we are dealing with  are positive multiple of each other.

\smallskip
    \noindent
     (1.6) The formula in either Theorem 1 or in Theorem 2 gives the Harish-Chandra parameter of $\mu \in Spec_H(\pi_\lambda^G)$ in term of the decomposition $\u = \u \cap \t_1(\Psi_\lambda) + \u \cap \t_2(\Psi_\lambda).$ However, other decompositions of $\u$ are: $\u =\z_h + \u \cap  \h_{ss} =\z_\l + \u \cap \l_{ss}.$  We now analyze the relation among these decompositions of $\u$  and its consequence in order to compute the Harish-Chandra parameter of $\pi_\mu^H.$ To begin with, we study an illustrative example. $\g_\mathbb R = \mathfrak{su}(2,1).$  We consider  $\Psi_\lambda =\{\alpha, \beta, \alpha + \beta \}$ with  $\alpha \in \Phi_c, \beta \in \Phi_n.$ $\h_\mathbb R = \t_\mathbb R + \mathfrak{sl}_2(\beta).$ Hence $\h_0 =\t_\mathbb R +\mathfrak{sl}_2(\alpha + \beta).$ For this system $$\t_1 =\z_K= Ker(\alpha)=\mathbb C H_{2\beta +\alpha}, \, \,\t_2 =\mathbb C H_\alpha,\,\,\,\,\, \z_\h = Ker(\beta)=\mathbb C H_{2\alpha +\beta}, \,\, \t \cap \h_{ss} = \mathbb C H_\beta.  $$ Thus, $\u =  \t_1(\Psi_\lambda) + \u \cap \t_2(\Psi_\lambda) =\z_\h + \u \cap  \h_{ss} $ are distinct orthogonal decompositions of $\u =\z_\l.$ More precisely, $\z_\h$ is not equal   $\t_1$ and is not equal to $ \z_\l \cap \t_2. $  Hence, in order to explicit the Harish-Chandra parameter $\mu$ as the sum of a central character for $\h$  plus  a Harish Chandra parameter for $\h_{ss}, $ we must carry out a  change of coordinates. Explicitly, $\mu= a (2\beta +\alpha) +b \alpha = \frac{a+b}{2} (2\alpha +\beta) + \frac{3a-b}{2}\beta.$   We now show that this picture prevails for most of the holomorphic systems $\Psi_\lambda$ and we analyze what happens when $\Psi_\lambda$ is not a holomorphic system.

     \noindent
     (1.6-a) For a nonholomorphic system $\Psi_\lambda$, that is, $\k_1= \k_1(\Psi_\lambda)$ is a simple Lie algebra, then $ \u \cap \h_{ss}=\t_1(\Psi_\lambda) + \t_2 \cap \h_{ss}  , \,\, \z_\h \subset \u \cap \t_2(\Psi_\lambda) \, \text{and} \z_\l \subset \t_2.$ Therefore,  $\mu = \mu_{\z_\h} + \mu_{\h_{ss}}, \mu_{\z_\h} \in \z_\h, \mu_{\h_{ss}} \in \h_{ss}  $  is easily computed from the decomposition $\mu =\mu_1 +\mu_2, \mu_j \in \t_j$  because $\mu_{\z_\h}= \mu_2$ restricted to $ \z_\h.$ \\  To show the claim we notice that $\pi_\lambda^G$ being an admissible representation of $H$ forces $\k_1 \subset \l.$ Hence, $\k_1$ is a simple Lie ideal in   $\l_{ss}.$ Thus, $ \t_1 \subset \u \cap \h_{ss}.$ Now, the orthogonal to $\t_1$ in $\u$ is equal to $\u \cap \t_2$ and $\z_\h$ is orthogonal to $\u \cap \h_{ss},$ hence, we have shown the second part of the claim.

\noindent
     (1.6-b)
    $\Psi_\lambda$ is a holomorphic system and $dim \z_\l =1 $. Then,  $\z_\h =\{0\}$ and $\z_{\h_0} =\{0\},$
      $\z_\l =\z_\k=\t_1(\Psi_\lambda),$  and  $\u \cap \t_2 (\Psi_\lambda)=\u \cap \l_{ss} .$  This hypothesis holds for: III-3(a) $1<k<n-1$, 3(b) $1<k<n-1$), 5 $(k \not= \{2,2n-1\})$, 9, 11, 13, 14.  \\ The first equality follows by inspection in table III. Since $res_H(\pi_\lambda^G)$ is an admissible representation of $H$ and $\Psi_\lambda$ is a holomorphic system we have  $\k_1=\z_\k=\t_1 \subset \l,$ hence, the second equality holds,    the third equality follows from that both members are the orthogonal to $\z_\l =\t_1$ in $\u.$

\noindent
(1.6-c)    $\Psi_\lambda$ is a holomorphic system, $\dim \z_\l =2$, and  $ \z_\h=\{0\}$    Then, $\z_\l = \t_1 + \z_\l \cap \t_2$ and
  $\u \cap \l_{ss} \subset \u \cap \t_2.$ The admissibility hypothesis forces $\t_1=\z_\k $ is contained in $\l.$ Orthogonality gives the inclusion. The hypothesis holds for:  III-3a ($k=1$), 3b $(k=n-1)$, 5a $(k=2$), 5b $(k=2n-1$), 6b ($1<m<n-1$), 7b, 10b, 12a.

\noindent
 (1.6-d)   $\Psi_\lambda$ is a holomorphic system,  $\dim \z_\l =2, $  and $dim  \z_\h =1 . $   Then, $\u=\t.$  The orthogonal decompositions of $\z_\l = \t_1 + \z_\l \cap \t_2 = \z_\h + \z_l \cap \h_{ss}$ are different  except for the pairs 3a, 3b, 5a, 5b, 6a $n=2m$, 7a $n=2m,$ for these pairs, we always have $\z_\h=\t_2 \cap \z_\l.$ We always have    $\t \cap \l_{ss} \subset \t_2$ and the orthogonal decomposition $ \t \cap \h_{ss}=\t \cap \l_{ss} + \z_\l \cap \h_{ss}.$ Thus, in order to obtain the decomposition $\mu = \mu_{\z_\h} + \mu_{\h_{ss}} $ from $\mu =\mu_1 +\mu_2$ we write $\mu_{\h_{ss}}= \mu_{\t \cap \l_{ss}} + \mu_{\z_\l \cap \h_{ss}}$ and notice that the inclusion $\t \cap \l_{ss} \subset  \t_2$ gives $\mu_{\t \cap \l_{ss}}$ is equal to $\mu_2$ restricted to $\t \cap \l_{ss}.$ The components $\mu_{\z_\h}, \mu_{\z_\l \cap \h_{ss}}$ are computed from $\mu_1 $ and from $\mu_2$ restricted to $\z_\l \cap \t_2$ as in the example $(SU(2,1), TSL_2(\beta)).$ The hypothesis holds for: III-1a-1b ($k=m, 1\leq l <n$ or $ 1 \leq k <m, l=n$ ), 2a-2b ($k=n, 1\leq l <n$ or $ 1 \leq k <n, l=n$ ), 3a ($k=n-1$),  3b ($k=1$), 4a, 4b, 5a ($k=2n-1$), 5b ($k=2$), 6a ($0<m<n$), 6b ($m=1$ or $ m=n-1$), 7a ($0<m<n$), 8a, 8b, 10a, 12b.

\noindent
(1.6-e)   $\Psi_\lambda$ is a holomorphic system and $\dim \z_\l =3.$   The hypothesis holds for cases: III-1a ($1\leq k <m $ and $ 1\leq l <n$), 1b ($1\leq k <m $ and $ 1\leq n-l <n$),  2a ($1\leq k <n $ and $ 1\leq l <n$), 2b ($1\leq k <n $ and $ 1\leq n-l <n$), 3a ($n=2, k=1$)). Then $dim \z_\h =1$ and $\u=\t.$ We have i) $\z_\h$ is orthogonal to $\t_1$ if and only if III-1a $mk=nl$, or, III-1b $nk=m(n-l)$, or, III-2a $k=l,$ or, III-2b $ n=k+l.$ ii) $\z_\h$ is not contained in $\t_2,$ obviously $\z_\h$ is not contained in $\t_1.$    Then, the orthogonal decompositions of $\z_\l = \t_1 + \z_\l \cap \t_2 = \z_\h + \z_l \cap \h_{ss}$ are different. In case i) we have $\t_1 \subset \h_{ss} \cap \t,$ hence $\mu_{\z_\h}$ is equal to the restriction of $\mu_2$ to $\z_\h.$  In case ii) we have $\t \cap \l_{ss} \subset \t_2$,  as in previous cases, we write $\u= \u_{\z_\h}+ \mu_{\t \cap \h_{ss}} = \u_{\z_\l} + \u_{\l_{ss}}$ and $\u_{\z_\l}= \u_{\z_\h} + \u_{\z_\l \cap \h_{ss}}.$ Then, $\u_{\h_{ss}}$ restricted to $\t \cap \l_{ss}$ es equal to  $\u_{\l_{ss}}=\mu_2$ restricted to $\t \cap \l_{ss}.$ In order to compute $\mu_{\z_\h}$ and $\u_{\z_\l \cap \h_{ss}}$ we must carry out a computation similar to the case $(\mathfrak{su}(2,1), \t +\mathfrak{sl}_2(\beta)).$ The case III-3a ($n=2, k=1$) is the case III-2a $(n=2,k=1,l=1).$  The first two claims follows from inspection to table III. To verify i) and ii) we recall  notation in $\mathbf{I-1}.$ Then, $\h=\mathfrak{su}(k,l)+\mathfrak{su}(n-k,m-l)+\u(1), \h_0=\mathfrak{su}(k,n-l)+\mathfrak{su}(n-k,l)+\u(1)$ and
  \begin{multline*}
     \Phi(\mathfrak h, \mathfrak t)=\{ \pm (\epsilon_r -\epsilon_s), \pm(\delta_p -\delta_q), \pm (\epsilon_i -\delta_j),
    1 \leq r < s \leq k \,\mbox{or} \, k+1 \leq r < s \leq m,\, \\  1 \leq p <  q \leq l  \,  \mbox{or} \,\, l+1 \leq p <  q \leq n,\\  1\leq i \leq k\,\, \mbox{and}  \,\, 1 \leq j \leq l\, \mbox{or}\, k+1\leq i \leq m \, \mbox{and} \,  l+1 \leq j \leq n \}.
    \end{multline*}
Let $e_j, d_k$ denote the elements of the dual basis to $\epsilon_j, \delta_k.$ Then, $$\z_\h = \mathbb C [ (m-k+n-l)(\sum_{i\leq k} e_i +\sum_{j \leq l} d_j)-(k+l)(\sum_{i\geq k+1} e_i +\sum_{j \geq l+1} d_j)],$$  $$ \t_1(\Psi_\lambda)=\z_\k =\mathbb C [ n(\sum_{i\leq m} e_i)- m(\sum_{j \leq n} d_j)]$$   and the inner product of the pointed out generators for $\z_\h, \z_\k$   is equal to $(nk-ml)(m+n).$ Thus, i) and ii) follow.

\section{More on the root structure for $(G,H)$}

 As in the previous sections, we assume $(G,H)$ is a symmetric pair and $\u=\h \cap \t$ is a Cartan subalgebra for $\h.$ We also assume $\pi_\lambda^G$ has an admissible restriction to $H.$ We  derive  consequences of the hypothesis and set up some notation. \\ For a linear subspace  $\s $ of $\g$  so that $[\u, \s] \subset \s $ or a complex vector space $\s$ where $\u$ acts by semisimple linear operators,  we denote by $\Phi(\mathfrak s, \u)$ the multiset of weights of $\u$ in $\s.$ That is, $\Phi(\mathfrak s, \u)$ is the set of weights of $\u$ in $\s,$ each counted with  multiplicity. Owing to our hypothesis we have that $\sigma (\Psi_\lambda)=\Psi_\lambda, $  $\k_1 \subset \l$ and $\sigma (\Delta) =\Delta, $ hence,  we may choose a system of positive roots $\Delta (\k /\l, \u)=\Delta(\k_2 /(\l \cap \k_2), \u)$ for $\Phi (\k /\l, \u)= \Phi (\k_2 /(\l \cap k_2), \u)$ naturally associated to $\Delta.$  For  $w \in W_K,$ as in \cite{dv},  we define the multiset  $$ S_w^H := [q_\u (w (\Psi_\lambda \cap \Phi_n)) \cup q_\u( \Delta \backslash \Phi_\mathfrak z)]\backslash \Phi(\h, \u)  .$$
Here, $\Phi_\z $ is the set of roots that vanishes on $\u.$  However, since we assume $\t$ is a maximally compact Cartan subalgebra for the symmetric pair $(\k ,\l)$ we have that $\Phi_\z $ is the empty set. The admissibility hypothesis yields the decompositions    $ L=K_1 (L\cap K_2), $  and  $\k=\k_1 \oplus \k_2. $ These decompositions give rise to the decomposition $W_{K}=W_{K_1} \times W_{K_2}.$
\begin{lem}  For $ W_K \ni w=ts, t \in W_{K_1}, s \in W_{K_2},$ we have the equality of the multisets
$$S_w^H = t [\Psi_{H_0, \lambda} \cap \Phi_n(\h_0,\u)] \cup \Delta(\k_2/ (\l \cap \k_2), \u).$$

\end{lem}

{\it Proof:}     When  $U=T$ we have defined $\Psi_{H_0,\lambda} :=\Psi_\lambda \cap \Phi(\h_0, \t)$ and $\Psi_{H,\lambda} :=\Psi_\lambda \cap \Phi(\h, \t).$ We have the disjoint union $\Phi_n= \Phi(\h, \t)_n \cup \Phi(\h_0, \t)_n = \Phi_n(\h, \t) \cup \Phi_n(\h_0, \t).$

 \begin{alignat*}{2}
    S_w^H &= (ts (\Psi_\lambda \cap \Phi_n)) \backslash \Phi(\h, \t) \cup \Delta(\k/\l, \t) &&\qquad\text{} \\
    &= (t (\Psi_\lambda \cap \Phi_n)) \backslash \Phi(\h, \t) \cup \Delta(\k/\l, \t) &&\qquad\text{by (1.4)} \\
    &=[t ((\Psi_\lambda \cap \Phi(\h,\t)_n) \cup (\Psi_\lambda \cap \Phi(\h_0,\t)_n))] \backslash \Phi(\h, \t)  &&\qquad\text{} \\ & \phantom{xxxxxxxxxxxxxxxxxxxxxxxxxxxx}\cup \Delta(\k_2/\l\cap \k_2 , \t) &&\qquad\text{ $\k_1 \subset \l$}\\
&=t [((\Psi_\lambda \cap \Phi(\h,\t)_n) \cup (\Psi_\lambda \cap \Phi(\h_0,\t)_n)) \backslash \Phi(\h, \t) ] &&\qquad\text{$t\Phi(\h) \subset \Phi(\h)$} \\ & \phantom{xxxxxxxxxxxxxxxxxxxxxxxxxxxx}\cup \Delta(\k_2/ (\l \cap \k_2), \t) &&\qquad\text{}\\
    & = t [\Psi_\lambda  \cap \Phi_n (\h_0, \t)] \cup \Delta(\k_2 / (\l \cap \k_2), \t). &&\qquad\text{}
 \end{alignat*}

 and     we conclude the proof of the lemma for the case $U=T.$

In order to show the   lemma when $U\not= T$ we  explicit some more  structure.

\noindent
(2.1) To begin with we note    the multiset $q_\u(\Phi(\g, \t)) $ is equal to  the multiset $\Phi(\h,\u) \cup \Phi(\h_0, \u) \cup \Phi(\frak q \cap \k, \u)$ and the equality of multiset $q_\u (\Phi_n)=\Phi_n (\h, \u) \cup \Phi_n(\h_0, \u) .$ Because $\sigma (\Psi_\lambda)=\Psi_\lambda$ we have $q_\u (\alpha) \in \Psi_{H,\lambda} \cup \Psi_{H_0, \lambda}$ for $\alpha \in \Psi_\lambda .$ Hence, the above equalities of multisets are also true when we replace everywhere $\Phi$ by $\Psi.$  Since $(G,H)$ is a symmetric pair in \cite{vothesis}, page 6 it is shown: for $\gamma_i \in \Phi(\g, \t), $ $q_\u(\gamma_1)=q_\u (\gamma_2)$ if and only if $\gamma_1=\gamma_2$ or $\gamma_1 =\sigma \gamma_2.$ Hence, $q_\u(\alpha) $ has multiplicity one in $\Phi(\g, \u)$ if and only if the root space of $\alpha$ is stable under $\sigma.$ Thus, as sets we have $\Psi_n(\h,\u)\cap \Psi_n (\h_0, \u)=\{q_\u(\beta), \beta \in \Psi_n \, \text{and} \, \beta \not= \sigma \beta \},$ and as multisets we have the equality $q_\u(\Psi_\lambda \cap \Phi_n) \backslash \Phi(\h, \u) = (\Psi_{H_0, \lambda})_n.$ Moreover, because of (1.4), for $s \in W_{K_2}, s(\Psi_\lambda \cap \Phi_n) =\Psi_\lambda \cap \Phi_n $ and since $K_1 \subset L$  we have $t\Phi(\h,\u)= \Phi(\h,\u)$ for $t \in W_{K_1}.$ Thus, $q_\u( ts (\Psi_\lambda \cap \Phi_n))\backslash \Phi(\h, \u)= t (\Psi_{H_0, \lambda})_n.$  The equalities we have obtained in the previous paragraphs justify the steps in:

\begin{equation*}
\begin{split}
S_w^H & =[t[ q_\u(\Psi_\lambda \cap \Phi_n)  ] \cup q_\u (\Delta)]\backslash \Phi(\h, \u) \\ &= t[ q_\u(\Psi_\lambda \cap \Phi_n) \backslash \Phi(\h, \u)] \cup \Delta (\k_2/\l \cap \k_2, \u) \\ &= t (\Psi_{H_0, \lambda})_n \cup \Delta (\k_2/\l \cap \k_2, \u).
\end{split}
\end{equation*}
This completes the proof of the lemma.
\begin{flushright}
$\Box$
\end{flushright}

 \noindent
 (2.1b) For $s  \in  W_{L \cap K_2} $ we have $s (\Psi_{H_0,\lambda} \cap \Phi_n(\h_0,\u))= \Psi_{H_0, \lambda} \cap \Phi_n (\h_0, \u).$

When $U=T,$ the equality readily follows from (1.4). When,  $U\not= T$ it can be done by direct computation. However, it also  follows from (1.4) and the following facts shown in  \cite{vothesis}, page 6.  For a root $\alpha \in \Phi(\g, \t) $ (resp. $\alpha \in \Phi(\h, \u)) $ let $S_\alpha^G $ (resp. $S_\alpha^H$) denote the reflection about $\alpha$ in $\t^\star.$ (resp. in $\u^\star).$

\noindent
(2.2) Let $\alpha $ be a root for $(\g, \t)$ and assume $q_\u(\alpha) \in \Phi(\h,\u).$  If $\alpha + \sigma(\alpha)$ is not a root for $(\g, \t).$  Then, for every $\gamma \in \t^\star,$   we have $$ S_{q_\u(\alpha)}^H (q_\u(\gamma)) = q_\u (S_{\alpha}^G S_{\sigma(\alpha)}^G (\gamma)).$$

\noindent
(2.3)   When $\alpha$ is a root in $(\g, \t)$ so that $\alpha =\sigma (\alpha)$ or $\alpha + \sigma (\alpha)$ is a root in $\g.$ We have the obvious equality $S_\alpha^H (q_\u(\gamma)) = q_\u ( S_\alpha^G (\gamma)).$ Now, 2.1b follows.

 \section{Proof of Theorem 1}

 Let $\pi_\lambda^G$ be a square integrable irreducible representation   which has an admissible restriction to $H.$ We want to show theorem 1 for this representation. To begin with,  we write several formulae  to compute the multiplicity functions $m^H(\lambda, \mu), \- m^{H_o,L}(\mu, \nu),\- m^{K_2,L \cap K_2}(\xi, \nu_2) .$ In order to write down the formulae we recall notation from \cite{dhv},\cite{dv}.

For  $ \nu \in i\mathfrak u_\mathbb R^\star, \, $ let $  \delta_\nu$ denote the Dirac measure attached to $\nu$ and for $\nu \not= 0,$ we consider the Heaviside discrete measure $y_\nu :=\sum_{n\geq 0} \delta_{\frac{\nu}{2} +n\nu}.$ For a strict multiset $S=\{\sigma_1, \dots, \sigma_n \} \subset \u^\star $ we define $y_S = y_{\sigma_1} \star \dots \star y_{\sigma_n}.$ Let $P_U$ denote the weight lattice for $U.$ Therefore, the set of Harish-Chandra parameters for $H$ (resp. $H_0$) is contained in $P_U +\rho_H$ (resp. $P_U +\rho_{H_0}$).  We would like to point out that either the set of Harish-Chandra parameters for $H$ (resp. $H_0$)  or $P_U +\rho_H$ (resp. $P_U +\rho_{H_0}$) are invariant under the action of the Weyl group of $L.$

Let $\mu $ be a Harish-Chandra parameter for $H_0 $  and dominant for $\Delta_0,$  for each $L-$type of $\pi_\mu^{ H_0},$ let $\nu \in i\mathfrak u_\mathbb R^\star $ denote the representative of its infinitesimal character which is dominant with respect to $\Delta_0$ and let $m^{ H_0, L}(\mu, \nu)$ denote the multiplicity of $\pi_\mu^L$ in $\pi_\lambda^G.$ We extend $m^{ H_0, L}(\mu, \nu)$ to $P_u +\rho_H$ by the rule $m^{ H_0, L}(\mu, w\nu)=\epsilon (w)m^{ H_0, L}(\mu, \nu)$ for $ w \in W_L.$ We denote by $\Psi_{ H_0, \mu} $ the system of positive roots in $\Phi(\h_0, \u)$ determined by $\mu.$ Finally, let  $\gamma_1, \dots,\gamma_r$ denote an enumeration of the noncompact  roots in $\Psi_{H_0, \mu}.$ Then, Duflo-Heckman-Vergne \cite{dhv} have shown,

$$ \sum_{\nu \in i\u_\mathbb R^\star } m^{ H_0, L}(\mu, \nu) \delta_\nu = \sum_{w \in W_L } \epsilon(w) \delta_{w\mu} \star y_{w\gamma_1}\star \cdots \star y_{w\gamma_r} \eqno(dhv) $$
The above series converges absolutely in the space of distributions of Schwartz for $i\u^\star.$
\noindent
The admissibility hypothesis on $res_H(\pi_\lambda^G)$  implies $K_1 \subset L $ (cf.(0.3)) and hence,  forces the decompositions  $L=K_1 (L\cap K_2)$ and the decomposition $\u =\t_1 +\u \cap  \t_2, \u^\star \ni \mu=\mu_1 +\mu_2, \mu_1 \in \t_1^\star, \mu_2 \in \t_2^\star. $  Thus, the right hand side of  formula (dhv)   is equal to $$ \sum_{t \in W_{K_1}, w_2 \in W_{L\cap K_2}} \epsilon (t) \epsilon (w_2) \delta_{t \mu_1} \star \delta_{w_2 \mu_2} \star y_{tw_2 (\Psi_{H_0, \mu})_n} .$$
We apply the above formula to Harish-Chandra parameters $\mu=(\lambda_1, \nu_2)$ as the one considered in  (1.5).  Hence, (2.1b) further simplifies the right hand side of (dhv) to $$ \sum_{t \in W_{K_1}, w_2 \in W_{L\cap K_2}} \epsilon (t) \epsilon (w_2) \delta_{t \mu_1} \star \delta_{w_2 \mu_2} \star y_{t (\Psi_{H_0, \lambda})_n} . \eqno(dhvs)$$
\noindent

In (1.4b) we noticed $\lambda_2$ is a Harish-Chandra parameter for $K_2$, hence,   for $res_{L\cap K_2}(\pi_{\lambda_2}^{K_2}).$ we have the identity
 \begin{multline*}\sum_{\nu_2 \in Spec_{L\cap K_2} (\pi_{\lambda_2}^{K_2})} m^{K_2, L\cap K_2}(\lambda_2, \nu_2) \sum_{s \in W_{L\cap K_2}} \epsilon (s) \delta_{s \nu_2} \\ = \epsilon_{12} \ \sum_{w_2 \in W_{K_2}} \epsilon(w_2) \delta_{q_\u(w_2 \lambda_2)} \bigstar y_{\Delta(\k_2 / (\l  \, \cap \k_2), \u)}.
 \end{multline*}
Here, $\epsilon_{12} =\epsilon(\Delta \cap \Phi(\k_2, \t), \Delta( \k_2 / \l \cap \k_2, \u)).$

  The hypothesis $\pi_\lambda^G$ is an admissible representation when restricted to  $H$ we give us  the equality $$ res_H(\pi_\lambda^G) = \sum_{\mu \in Spec_H(\pi_\lambda^G)}m^H(\lambda, \mu) \, \pi_\mu^H.$$ We notice in \cite{oshi} it is shown that any Harish-Chandra parameter in $Spec_H(\pi_\lambda^G)$ is  dominant for the system  $\Psi_{H,\lambda}.$ We extend $m^H(\lambda, \mu)$ to $P_U + \rho_H$ to be a skew symmetric function, that is, $m^H(\lambda, w\mu) =\epsilon(w) m^H(\lambda, \mu) $ for $w \in W_L.$  As before, for $ w \in W_K, $ let $$S_w^H = [q_\u (w (\Psi_\lambda \cap \Phi_n)) \cup q_\u( \Delta \backslash \Phi_\mathfrak z)]\backslash \Phi(\h, \u)= q_\u (w (\Psi_\lambda \cap \Phi_n)) \backslash \Phi(\h, \u) \cup \Delta(\k/\l, \u) $$ Then, in \cite{dv} is proved
$$ \sum_{\mu \in i\u_\mathbb R^\star } m^H(\lambda , \mu) \delta_\mu = \pm \sum_{ w \in   W_K} \epsilon(w) \, \delta_{q_\u(w\lambda)} \star y_{S_w^H} \eqno(rh) $$
The above series converge absolutely in the space of distributions of Schwartz for $i\u^\star.$
Actually, in (rh)  the summation on the right hand side  is computed over the group $W_\z /W_K$ where $W_\z$  is the Weyl group of the root system $\Phi_\mathfrak z $ of those roots in $\Phi(\g, \t)$ which vanishes on $\u.$ The hypothesis $(G,H)$ is a symmetric pair implies $(K,L)$ is a symmetric pair, hence  $W_\z $ is equal to the trivial group. Since $\t_1 $ is contained in $\l, $ and  $\lambda =(\lambda_1, \lambda_2), $ for every $w_1 \in W_{\k_1}, w_2 \in W_{\k_2},$ we have  $$ q_\u(w_1 w_2 \lambda)= w_1 \lambda_1+ q_\u(w_2\lambda_2).$$

\noindent
Also, for $\nu_2 \in  \t_2^\star \cap \u^\star, t \in W_{K_1}, s \in W_{K_2}$ we have  $ts(\lambda_1 +\nu_2)= t \lambda_1 + s \nu_2 \, \text{and} \, \delta_{t \lambda_1} \star \delta_{s \nu_2} = \delta_{t s (\lambda_1 +\nu_2)}.$
 The previous equalities, Lemma 1, (2.1b), (dhvs) and the considerations in the previous paragraphs,  justify the following  transformations on  the right hand side in $(rh)$
\begin{multline*}
 \sum_{\mu \in Spec_H(\pi_\lambda^H)} m^H(\lambda, \mu) \, \delta_\mu\\
=\sum_{t \in W_{K_1}} \epsilon(t) \delta_{t \lambda_1} \bigstar \sum_{w_2 \in W_{K_2}} \epsilon (w_2) \delta_{q_\u (w_2 \lambda_2)} \bigstar y_{\Delta(\k_2/ (\l  \,\cap \k_2), \u)} \bigstar y_{t (\Psi_{H_0, \lambda})_n} \\ =\sum_{t \in W_{K_1}} \epsilon(t) \delta_{t \lambda_1}  \bigstar y_{t (\Psi_{H_0, \lambda})_n} \bigstar \sum_{\nu_2 \in Spec_{L\cap K_2}(\pi_{\lambda_2}^{K_2}), \, s \in W_{L\cap K_2}} m^{K_2, L\cap K_2} (\lambda_2, \nu_2) \epsilon(s) \delta_{s \nu_2} \\ = \sum_{\nu_2 } m^{K_2, L\cap K_2} (\lambda_2, \nu_2) \sum_{t \in W_{K_1}, s \in W_{ L \cap K_2}} \epsilon(t) e(s) \delta_{t  s \lambda_1} \bigstar  \delta_{t s \nu_2} \bigstar y_{t  s (\Psi_{H_0, \lambda})_n}
\\= \sum_{\nu_2 \in Spec_{L\cap K_2}(\pi_{\lambda_2}^{K_2}) } \,\, \sum_{\xi \in i\u_\mathbb R^\star}  m^{K_2, L\cap K_2} (\lambda_2, \nu_2)\, m^{H_0, L}((\lambda_1, \nu_2), \xi)  \,\, \delta_{\xi}
\end{multline*}

Since $P_U$ is a discrete subset of $\u^\star, $ and the above series converge absolutely in the topology of the space of distributions on $i\u^\star,$   we have  shown Theorem 1.

As a consequence, we obtain particular cases of a more general result shown in \cite{kobhol}.
\begin{pro} Assume $\Psi_\lambda$ is a holomorphic system. Then, whenever $res_H(\pi_\lambda^G)$ is an admissible representation of $H,$  there exists a constant $C<\infty$ so that  $m^H(\lambda, \mu))\leq C $ for every $\mu \in Spec_H(\pi_\lambda^G.)$
\end{pro}
{\it Proof:} It follows from the hypothesis  $\Psi_\lambda$ is a holomorphic system that $\Psi_{H_0,\lambda}$ is also a holomorphic system. Hence, the subspaces $$\p_0^-=\sum_{-\beta \in (\Psi_{H_0,\lambda})_n} \h_\beta,\,\,\,\,\,\, \p_0^+=\sum_{\beta \in (\Psi_{H_0,\lambda})_n} \h_\beta $$ are abelian subalgebras and we have the $L-$invariant and  direct sum decomposition $\p \cap \p  =\p_0^+ + \p_0^-.$  Let $(\pi_{\mu +\rho_n^0}^L , W_{\mu +\rho_n^0}) $ denote the lowest $L-$type of the Discrete series $\pi_\mu^{H_0}.$ For a vector space $V$, the symmetric algebra for $V$ is denoted by $S(V).$  For a Harish-Chandra parameter for $\h_0$  and dominant for $\Psi_{H_0,\lambda},$  a result of Harish-Chandra \cite{knappov}   gives us  that $res_L(\pi_\mu^{H_0}) $ is equivalent to  $S(\p_0^+ ) \otimes W_{\mu +\rho_n^0}. $ Since each $L-$irreducible factor of $S(\p_0^+ )$ has  multiplicity one, we have, owing to a tensor product argument, that there exist a constant $C<\infty$ so that $m^{H_0, L}(\mu, \xi)\leq C$ for every Harish-Chandra parameter $ \xi$ for $L.$ The formula in Theorem 2 concludes the proof of the proposition.

\begin{cor} For a scalar holomorphic discrete series, then $m^H(\lambda, \mu)=1$ for every $\mu \in Spec_H(\pi_\lambda^G).$

\end{cor}
\section{Analysis of structure of $L-$types}

       As in the previous setting  we assume $\pi_\lambda^G$ has an admissible restriction to $H$  as well as that  $(G,H)$ is a symmetric pair. We also assume that the lowest $K-$type is an irreducible representation of $K_1.$ Then, (4.1) yields that the Harish-Chandra parameter $\lambda=(\lambda_1, \rho_{K_2})$ and the lowest $K-$type is $(\pi_{(\lambda_1+ \rho_n^\lambda, \rho_{K_2} )}^K, V_{(\lambda_1+ \rho_n^\lambda, \rho_{K_2} )}^K).$ Theorem 2 shows that $Spec_H(\pi_\lambda^G)=Spec_{L}(\pi_{(\lambda_1, \rho_{K_2 \cap L})}^{H_0})$ and $m^H(\lambda, \mu)=m^{H_0, L}( (\lambda_1, \rho_{K_2 \cap L}), \mu).$ Now, (4.1) together with the fact that the lowest $L-$type of a discrete series representation for $H$ has multiplicity one yield \begin{multline*}  dim Hom_H(\pi_{(\lambda_1 +\rho_0, \rho_{L\cap K_2})}^H , res_H(\pi_\lambda^G))=m^H(\lambda, (\lambda_1 +\rho_0, \rho_{L\cap K_2})) \\ =m^{H_0,L}((\lambda_1, \rho_{K_2 \cap L}), (\lambda_1 +\rho_0, \rho_{K_2 \cap L})) =1    \end{multline*}  Here, and from now on, in order to avoid cumbersome notation, we write  $\rho_0$(resp $\rho_1)$  for one half of the sum of the noncompact in $\Psi_{H_0,\lambda}$(resp. $\Psi_{H,\lambda}).$ Thus, when $U=T$ we have $\rho_n^\lambda=\rho_1 +\rho_0.$ We now show
       \begin{pro} The lowest $L-$type $\pi_{(\lambda_1 +\rho_n^\lambda, \rho_{L\cap K_2})}^L $ of $\pi_{(\lambda_1 +\rho_0, \rho_{L\cap K_2})}^{H} $ has multiplicity one in $res_L(\pi_\lambda^G).$
       \end{pro}
       {\it Proof:} Let $\pi_\mu^H$ be an irreducible subrepresentation of $res_H(\pi_\lambda^G)$ which contains a copy of the $L-$type $\pi_{(\lambda_1 +\rho_n^\lambda, \rho_{L\cap K_2})}^L .$  Because of (4.1) we have $$ (\lambda_1+\rho_n^\lambda, \rho_{L\cap K_2})= \mu +\rho_1 +D_1,$$ here $D_1$ stands for a sum of noncompact roots in $\Psi_{H, \lambda}.$ \\ Theorem 2 let us write $$\mu=(\lambda_1 +\rho_0, \rho_{L\cap K_2})+ B_0$$ where $B_0$ is a sum of noncompact  roots in $\Psi_{H_0, \lambda}.$   Hence, $D_1+B_0=0.$ The fact that $\Psi_\lambda$ is invariant under $\sigma$ forces $D_1=B_0=0$ and  proposition 2 follows.

\medskip
\noindent
For the next lemma we keep the assumption that the lowest $K-$type of $\pi_\lambda^G$ is an irreducible representation $(\pi_{(\lambda_1 +\rho_n^\lambda, \rho_{K_2})}^K, V_{(\lambda_1 +\rho_n^\lambda, \rho_{K_2})}^K)$ for $K_1.$
\begin{lem} $\mathcal U(\h_0)V_{(\lambda_1 +\rho_n^\lambda, \rho_{K_2})}^K$ is an irreducible $(\h_0,L)-$submodule and has multiplicity one in $res_{H_0}(\pi_\lambda^G).$
\end{lem}
{\it Proof:} Since we are dealing with restriction of discrete series representation, \cite{dv} we have that the underlying Harish-Chandra  module for $\pi_\lambda^G$ is an admissible representation of $L.$ Thus, $\mathcal U(\h_0)V_{(\lambda_1 +\rho_n^\lambda, \rho_{K_2})}^K$ is an admissible $L-$module. Besides, since the Harish-Chandra module for $\pi_\lambda^G$ is admissible as representation of $\h_0$ we have that every vector en the Harish-Chandra module for $\pi_\lambda^G$  is $\z_{\mathcal U(\h_0)}-$finite. Thus, we are in the hypothesis of Theorem 4.2.1 in \cite{wallb1}. Therefore, $\mathcal U(\h_0)V_{(\lambda_1 +\rho_n^\lambda, \rho_{K_2})}^K$ is an $\h_0-$module of finite length. Whence, the closure of $\mathcal U(\h_0)V_{(\lambda_1 +\rho_n^\lambda, \rho_{K_2})}^K$ is equal to a finite orthogonal sum $V_1 +\dots +V_s$ of discrete series representations $V_j$ for $H_0.$ Since, $V_{(\lambda_1 +\rho_n^\lambda, \rho_{K_2})}^K$ is a cyclic generator for $\mathcal U(\h_0)V_{(\lambda_1 +\rho_n^\lambda, \rho_{K_2})}^K$ the orthogonal projector onto $V_i$ takes on nonzero values on the $L-$irreducible subspace $V_{(\lambda_1 +\rho_n^\lambda, \rho_{K_2})}^K. $ Proposition 2 forces $s=1$ and  lemma 2 follows.

\medskip

As previously, we consider $\pi_\lambda^G$ to be an admissible representation for $H.$  Let $$\mu_1, \mu_2, \dots    $$ denote the Harish-Chandra parameters of the distinct  irreducible $H-$factors of $res_H(\pi_\lambda^G), $ we know \cite{gw}, \cite{oshi} every $\mu_j$ is dominant with respect to $\Psi_{H,\lambda}.$   \\ For each $\pi_{\mu_j}^H, j=1, 2, \cdots $     the lowest $L-$type in the sense of Schmid-Vogan of $\pi_{\mu_j}^H$ is $\pi_{\mu_j + \rho_1}^L.$ The next proposition holds for the pairs  $ (\mathfrak{su}(m,n), \mathfrak{su} (m,l) +\mathfrak{su} (n-l)+\u(1)),$ $(\mathfrak{so}(2,2n), \u(1,n)),$ $(\mathfrak{so}^\star(2n),\u(1,n-1)),$  $(\mathfrak e_{6(-14)}, \mathfrak{so}(2,8)+\mathfrak{so}(2)).$

\begin{pro} We further assume $\Psi_\lambda$ is a holomorphic system. Then, for $r \not= j,$    $$ Hom_L(\pi_{\mu_j +\rho_1}^L, res_L(\pi_{\mu_r}^H))=\{0 \}$$
\end{pro}

{\it Proof:}    For each of the  pair   $ (\mathfrak{su}(m,n), \mathfrak{su} (m,l) +\mathfrak{su} ( n-l)+\u(1)),$ $(\mathfrak{so}(2,2n), \u(1,n)),$ $(\mathfrak{so}^\star(2n),\u(1,n-1)),$ $(\mathfrak e_{6(-14)}, \mathfrak{so}(2,8)+\mathfrak{so}(2)),$ we show

   \begin{lem}
There exists $x \in i\t, y \in i\t$ such that \begin{center}  for $ \alpha \in  \Psi_\lambda \cap \Phi_n (\h, \t), \,\,  \,\, \alpha(x)>0 \,\, \text{and} \,\, \alpha (y)=0,$

for $ \alpha \in  \Psi_\lambda  \cap \Phi_n (\h_0, \t), \,\, \,\, \alpha(x)=0 \,\, \text{and} \,\, \alpha (y)>0.$
\end{center}
 \end{lem}
  We note that the existence of a nonzero $x \in i\t$ so that $\alpha(x)=0 $ for every noncompact root in $\Phi(\h,\t)$  forces  the center of $\h$ to have positive dimension because of the equality    $[\p, \p]=\k.$ However, it readily follows that the lemma is not true  for the pair $(\mathfrak{su}(m,n), \mathfrak{su} (k,l) +\mathfrak{su} (m-k, n-l)+\u(1))   $ and its dual when $1 \leq k <m, 1 \leq l <n.$ A proof of  lemma 3 is done case  by case after we verify proposition 3.

\smallskip
\noindent
(4.1) We recall the following result proved in \cite{schannals}.
The highest weight of any  $K-$type of $\pi_\lambda^G$ is equal to $\lambda +\rho_n^\lambda -\rho_c +\beta_1 + \dots \beta_s, $ where $\beta_i, i=1, \dots, s$ are noncompact roots in $\Psi_\lambda.$\\

\smallskip
\noindent
(4.2) We notice that whenever the lowest $K-$type of $\pi_\lambda^G$ is an irreducible representation of $K_1,$  the Harish-Chandra parameter   is  $\lambda=(\lambda_1, \rho_{K_2}) $ with $\lambda_1 +\rho_n^\lambda $ a Harish-Chandra parameter for $K_1.$  This is so, because from (1.4a) we have that $\rho_n^\lambda$ lies in $\t_1^\star, $  also $\k=\k_1 +\k_2$ a direct sum of ideals. Thus (4.2) follows.\\
    We now show proposition 3,  a  Harish-Chandra parameters for a $L-$ of types of $\pi_{(\lambda_1, \rho_{L\cap K_2})}^{H_0}$ is of the shape $$(\lambda_1 , \rho_{L\cap K_2})+\rho_0+ B$$ where $B$ is a sum of roots in $\Psi_{H_0,\lambda} \cap \Phi_n (\h_0, \t).$ Hence, we may and will  order the Harish-Chandra parameters of the $L-$types of  $\pi_{(\lambda_1, \rho_{L\cap K_2})}^{H_0}$ in an increasing way according to the value each of them takes on $y.$ Because, of Theorem 2, this gives an order $$\mu_1 < \mu_2 < \dots . $$   Also, the Harish-Chandra parameter of an  $L-$type of $\pi_{\mu_r}^H$ is $$\mu_r + \rho_1 + C$$ where $C$ is a sum of noncompact roots in $\Psi_{H,\lambda}.$       Thus, if $Hom_L(\pi_{\mu_j +\rho_1}^L , res_L(\pi_{\mu_r}^H))$ is nontrivial, we have the equality $ \mu_j +\rho_1 = \mu_r +\rho_1 + C $ where $C$ is a sum of roots in $\Psi_{H,\lambda} \cap \Phi_n.$ Since $C(y)=0,$ we obtain $r=j$  and we have shown   proposition 3.

\smallskip
\noindent
In order to justify lemma 2 for the case  of  $(\mathfrak{su}(m,n), \mathfrak{su} (m,l) +\mathfrak{su} (n-l)+\u(1))$ we recall  notation in $\mathbf{I-1},$ we set $e_j, d_s$ to be  the dual basis to $\epsilon_j, \delta_s.$  Then, $x=\sum_1^m e_j  -\sum_1^l  d_j +\sum_{l+1}^n  d_j$ and  $y=\sum_1^m e_j  +\sum_1^l  d_j -\sum_{l+1}^n  d_j$  verify the claim for this pair.

\noindent
 $\mathbf{III-4}  \, (\mathfrak{so}(2m,2), \u(m,1)).  $ We refer to  notation  in $\mathbf{II-5}.$ Then $\Phi (\h, \t)=\{\pm (\epsilon_j -\delta_1), \pm(\epsilon_s -\epsilon_k), s\not=k  \},$  $ \Phi (\h_0, \t)=\{\pm (\epsilon_j +\delta_1), \pm(\epsilon_s -\epsilon_k), k\not=s  \}  .$ For the holomorphic system  $\Psi_+ , $  $x=-\sum_{1}^m  e_j +d_1, y= \sum_{1}^m  e_j +d_1$ verify the claim.

\noindent
$(\mathfrak{so}^\star(2n),\u(1,n-1)).$ Notation as in $\mathbf{III-6}.$ $x=e_1, y=-e_1 +(e_2 +\dots +e_n).$  For  $(\mathfrak{so}^\star(2n),\u(n-1,1))$ we do a similar choice.

\noindent
$\mathbf{III-8} \, (\mathfrak e_{6(-14)}, \mathfrak{so}(2,8)+\mathfrak{so}(2)).$  The Vogan diagram for a holomorphic system for $\mathfrak e_{6(-14)}$ is

\begin{figure}[h]
\begin{tikzpicture}
\draw (0.05,0) --(0.94,0);  \draw (1.04,0) --(1.97,0); \draw
(2.04,0) -- (2.97,0); \draw (3.03,0) -- (3.96,0) ;
\draw (2,0) -- (2,1);
\node[left] at (2,1) {$\alpha_2$};
\node at (2,1) {$\circ$};
\node at (2,0) {$\circ$};
\node at (0,0) {$\circ$};
\node at (1,0) {$\circ$};
\node at (3,0) {$\circ$};
\node at (4,0) {$\bullet$};
\node[below] at (0,0) {$\alpha_1$};\node[below] at (1,0) {$\alpha_3$};\node[below] at (2,0) {$\alpha_4$};\node[below] at (3,0) {$\alpha_5$}; \node[below] at (4,0) {$\alpha_6$};
\end{tikzpicture}
\end{figure}

\noindent
Then $\Phi_n(\mathfrak{so}(2,8)+\mathfrak{so}(2) ,\t)=\{ \sum_j a_j \alpha_j \in \Phi(\mathfrak e_6,\mathfrak t) :  a_1=0, a_6=1\}, \Phi_n(\h_0,  \t)=\{ \sum_j a_j \alpha_j \in \Phi(\mathfrak e_6,\mathfrak t) :  a_1=1, a_6=1\}$ and $\Phi(\l, \t)=\{ \sum_j a_j \alpha_j \in \Phi(\mathfrak e_6,\mathfrak t) :  a_1=a_6=0\}.$ Here, $\z_\h^\star=\mathbb C \Lambda_1, \z_{\h_0}^\star= \Lambda_1 -\Lambda_6.$ Hence $x=\Lambda_6 -\Lambda_1, y=\Lambda_1$ completes the proof of  lemma 2. Here, $\Lambda_j$ is the fundamental weight associated to $\alpha_j.$

\subsection{The subspace of lowest $L-$types in $res_H(\pi_\lambda^G)$}

We fix $(\pi_\lambda^G, V_\lambda^G)$ a discrete series representation, for which, we assume has an admissible restriction to $H.$ Since, we are dealing with discrete series representations and symmetric pairs    in  \cite{dv} we find a proof that $res_L(\pi_\lambda^G)$ is an admissible representation of $L$ and in \cite{kobic} a proof that the underlying Harish-Chandra module $\mathcal V_\lambda^G$ of $V_\lambda^G$ is an algebraic direct sum of irreducible  $(\mathfrak h, L)-$modules, namely  $$ \mathcal V_\lambda^G =\bigoplus_{\mu \in Spec_H(\pi_\lambda^G)} \,\, m^{H}(\lambda, \mu)\,\,\, \mathcal V_\mu^H. $$  Here, $ m^{H}(\lambda, \mu)\mathcal V_\mu^H$ denotes the isotopic component of $\mathcal V_\lambda^G$ associated to  the discrete series representation $\pi_\mu^H.$ For each subspace  $\mathcal V_\mu^H$ the subspace which affords the lowest $L-$type of $\pi_\mu^H,$ or equivalently, the subspace of vectors which behaves according to the representation  $\pi_{\mu +\rho_1}^L$ of $L$ is denoted by $\mathcal V_{\mu + \rho_1}^L.$ A  problem considered for \cite{os} is to analyze the structure of the subspace $$ \mathcal L_\lambda :=\bigoplus_{\mu \in Spec_H(\pi_\lambda)}\,\, m^{H}(\lambda, \mu)\,\,\, \mathcal V_{\mu+\rho_1}^L. \eqno(LL) $$ For a partial answer, we further assume  the lowest $K-$type of $\pi_\lambda^G $ restricted to $K_1$ is an irreducible representation.  We claim:

 \smallskip
 \noindent
 (4.3)  The respective representations of $L$ on the  subspaces   $\mathcal L_\lambda $ and on   $\mathcal U(\mathfrak h_0) \mathcal V_{(\lambda_1 +\rho_n^\lambda, \rho_{K_2})}^K $  are equivalent.

\noindent
Indeed, lemma 2 implies that  $\mathcal U(\mathfrak h_0) \mathcal V_{(\lambda_1 +\rho_n^\lambda, \rho_{K_2})}^K $ is the underlying Harish-Chandra module for $\pi_{(\lambda_1 +\rho_1 ,\rho_{L\cap K_2})}^{H_0},$  Theorem 2 concludes the proof of the claim.

\medskip
\noindent
To continue, we   assume $\pi_\lambda^G$ is a holomorphic representation.   We set $\p^+ =\sum_{\beta \in (\Psi_\lambda)_n} \g_\beta $ and $\p^- =\bar{\p^+}. $  Since $\z_\k$ is contained in $\l$  we also have $\h =\l + \h \cap \p^+ + \h \cap \p^-.$ Thus, the systems $\Psi_{H, \lambda}, \Psi_{H_0, \lambda}$ are holomorphic.  The representation of $K$ in  $$ (\mathcal V_\lambda^G)^{\mathfrak p^-} :=\{ v \in \mathcal V_\lambda^G : \pi_\lambda (Y)(v)=0, \,\,\forall \,\,Y \in \mathfrak p^- \}$$ is the realization of  $\mathcal V_{(\lambda_1 +\rho_n^\lambda, \rho_{K_2})}^K$ as a subspace of $\mathcal V_\lambda^G.$    Moreover,  the fact  $\Psi_{H,\lambda}$ is holomorphic   and that  for each constituent $\pi_\mu^H $ of $res_H(\pi_\lambda^G)  $ the Harish-Chandra parameter $\mu$ is  dominant for $\Psi_{H,\lambda}$  yield

$$ \mathcal V_\lambda^{\mathfrak h \cap \mathfrak p^-} = \mathcal L_\lambda.$$

\noindent
 Certainly $ \mathcal V_{(\lambda_1 +\rho_n^\lambda, \rho_{K_2})}^K \subset \mathcal V_\lambda^{\mathfrak h \cap  \mathfrak p^-} \cap \mathcal U (\mathfrak h_0) \mathcal V_{(\lambda_1 +\rho_n^\lambda, \rho_{K_2})}^K.$

 \begin{pro} We further assume $\pi_\lambda^G$ is scalar holomorphic discrete series representation. Then,

 i) If $[[\mathfrak h \cap \mathfrak p^-, \mathfrak h_0 \cap \mathfrak p^+], \mathfrak h_0 \cap \mathfrak p^+]=\{0\},$ then $\mathcal V_\lambda^{\mathfrak h \cap \mathfrak p^-} =\mathcal U(\mathfrak h_0) (\mathcal V_{(\lambda_1 +\rho_n^\lambda, \rho_{K_2})}^K).$

 ii)  If $[[\mathfrak h \cap \mathfrak p^-, \mathfrak h_0 \cap \mathfrak p^+], \mathfrak h_0 \cap \mathfrak p^+]\not= \{0\},$ then $\mathcal V_\lambda^{\mathfrak h \cap \mathfrak p^-} \not= \mathcal U(\mathfrak h_0) (\mathcal V_{(\lambda_1 +\rho_n^\lambda, \rho_{K_2})}^K).$

 \end{pro}
In a case by case checking, we verify that the unique pairs $(\g, \h)$ that satisfy  the hypothesis in i) are:\\
(4.4) $ (\mathfrak{su}(m,n), \mathfrak{su} (m,l) +\mathfrak{su} ( n-l)+\u(1)),$ $(\mathfrak{so}(2m,2), \u(m,1)),$  $(\mathfrak{so}^\star (2n), \u(1,n-1)),$   $(\mathfrak e_{6(-14)}, \mathfrak{so}(2,8)+\mathfrak{so}(2)).$ That is, the same pairs that satisfy Lemma 3. Whence, when consider case i)  $T$ is a Cartan subgroup of $H$.\\
 {\it Proof:} In order to verify i) we  show that left hand side of the equality is contained in $\mathcal V_\lambda^{\h \cap \p^-},$ in (4.3) we observed $\mathcal L_\lambda$ and $\mathcal U(\mathfrak h_0) (\mathcal V_{(\lambda_1 +\rho_n^\lambda, \rho_{K_2})}^K)$ have the same $L-$module structure, hence the reverse inclusion follows. In order to verify ii) we produce an element  in the left hand side which is not in $\mathcal V_\lambda^{\h \cap \p^-}.$  For each root $\alpha \in \Phi(\mathfrak g, \mathfrak t)$ we choose a nonzero root vector $Y_\alpha$ in the root space of $\alpha.$ Owing to our hypothesis we may write $\mathcal V_{(\lambda_1 +\rho_n^\lambda, \rho_{K_2})}^K=\mathbb C w$ with $w $ a nonzero vector. We denote by  $\Phi(\mathfrak p^+, \mathfrak t)=\{ \beta_1, \dots, \beta_q\}, $ where $\Phi(\h_0 \cap \p^+)=\{ \beta_1, \dots \beta_s \}.$ Since, $  \mathcal V_\lambda^G \equiv  U(\g)_{\mathcal U(\k +\p^-)} \otimes \mathcal V_{(\lambda_1 +\rho_n^\lambda, \rho_{K_2})}^K \equiv_K  S(\mathfrak p^+ )\otimes \mathcal V_{(\lambda_1 +\rho_n^\lambda, \rho_{K_2})}^K   $  the projection onto $\mathcal V_\lambda^G$ of the  set $$ \{ Y_{\beta_1}^{a_1} \cdots Y_{\beta_q}^{a_q} \otimes w , a_j \in \mathbb Z_{\geq 0}, j=1, \dots, q \} $$ is a  linear basis for $\mathcal V_\lambda^G.$   Thus, $\mathcal U(\mathfrak h_0) \mathcal V_{(\lambda_1 +\rho_n^\lambda, \rho_{K_2})}^K$ is equal to the subspace of $\mathcal V_\lambda^G$ spanned by the projection of $S(\mathfrak h_0 \cap \mathfrak p^+) \otimes \mathbb C w.$ Under the hypothesis in i) we verify that $$\pi_\lambda(Y_{-\beta_\h} ) [Y_{\beta_1}^{a_1} \cdots Y_{\beta_s}^{a_s} \otimes w]=0 \,\,\text{for\,all}\,\,Y_{-\beta_\h} \in \h \cap \p^-. $$ by induction on $a_1 +\dots +a_s.$ Here, $[D \otimes w]$ denotes the class of $D\otimes w \in \mathcal U(\g) \otimes V_{(\lambda_1 +\rho_n^\lambda, \rho_{K_2})}^K.$ \\  Since, $\pi_\lambda^G$ is an scalar holomorphic representation, we have $\pi_\lambda(Y_{-\beta_j})[1\otimes w]=0, \forall \,\,\beta_j \in \Psi_\lambda \cap \Phi_n $ and, $\pi_\lambda (Y_\gamma)[1\otimes w]=0 , $ for every $\gamma \in  \Phi_c.$ Hence, for $j \leq s,$ we have $\pi (Y_{-\beta_\h})[Y_{\beta_j} \otimes w]=[Y_{\beta_j}Y_{-\beta_\h} \otimes w + [Y_{-\beta_\h}, Y_{\beta_j}] \otimes w]=0$ because $ [Y_{-\beta_\h}, Y_{\beta_j}] \in \mathfrak q \cap \k $ which is contained in $\k_{ss}.$ In general, for $ a_1= \dots =a_{j-1}=0, a_j \geq 1, j\leq s$ we have $ \pi (Y_{-\beta_\h})[Y_{\beta_j}^{a_j} \dots Y_{\beta_s}^{a_s} \otimes w]=[Y_{\beta_j}Y_{-\beta_\h}Y_{\beta_j}^{a_j -1} \dots Y_{\beta_s}^{a_s} \otimes w + [Y_{-\beta_\h}, Y_{\beta_j}] Y_{\beta_j}^{a_j -1} \dots Y_{\beta_s}^{a_s} \otimes w]=0$ if $a_j=1, $ otherwise,  we  may assume $a_j \geq 2. $  In this  case,  the hypothesis in i)  states that $[Y_{-\beta_\h}, Y_{\beta_i}]$ commutes with $Y_{\beta_k}$ for $1\leq i \leq s, 1 \leq k \leq s.$ Hence, we have shown i). In order to show ii) when $U=T$  we compute  $ \pi (Y_{-\beta_\h})[Y_{\beta_j} Y_{\beta_r} \otimes w]=[Y_{\beta_j}Y_{-\beta_\h}Y_{\beta_r} \otimes w + [Y_{-\beta_\h}, Y_{\beta_j}]   Y_{\beta_r} \otimes w]=[ Y_{\beta_j}Y_{\beta_r}Y_{-\beta_\h}  \otimes w +Y_{\beta_j}[Y_{-\beta_\h},Y_{\beta_r}] \otimes w  + [[Y_{-\beta_\h}, Y_{\beta_j}]  , Y_{\beta_r}] \otimes w + Y_{\beta_r} [Y_{-\beta_\h}, Y_{\beta_j}]     \otimes w    ]= [[Y_{-\beta_\h}, Y_{\beta_j}]  , Y_{\beta_r}] \otimes w ].$ Now $ [[Y_{-\beta_\h}, Y_{\beta_j}]  , Y_{\beta_r}] \in \mathfrak h \cap \p^+ .$ Under the hypothesis in ii) we may choose $\beta_\h \in (\Psi_{H, \lambda})_n , \beta_j, \beta_r, \text{with} j,r \leq  s$ so that the triple bracket is nonzero, and we have shown ii) in case $U=T.$ When $U\not= T$ always the hypothesis in ii) holds because we check by the end of section 6 that there exists $\beta \in  \Phi_n(\h, \u) \cap \Phi_n(\h_0, \u) \not= \emptyset$  and  $X_{-\beta} \in \h_{-\beta}$ and $V_\beta \in (\h_0)_\beta $ so that  $[[X_{-\beta}, V_\beta], V_\beta] \not= 0.$   Whence,  $\mathcal U(\h_0) (1\otimes w)\not=\mathcal L_\lambda. $ Thus, we have shown proposition 4.

\section {Tensor product of holomorphic Discrete Series}

Let $G_0$ be a real simple Lie group so that the associated symmetric space is Hermitian. Let $K_0$ (resp. $T_0$)  denote a maximal compact subgroup (resp. maximal torus) of $G_0$ (resp. $K_0$). Thus, $T_0=Z_{K_0} T_0^s$ where $T_0^s$ is a maximal torus of the semisimple factor $K_0^s$ of $K_0.$  Let $\theta_0$ denote the Cartan involution of $G_0$ associated to $K_0.$  We set  $G:=G_0 \times G_0,  $ hence $K:=K_0 \times K_0$  (resp. $T:=T_0 \times T_0$) is maximal compact subgroup  (resp. is a compact Cartan subgroup) of $G.$ In $G$ we consider the involution $\sigma (x,y) = (y,x).$ Thus, the fix point subgroup $H$ of $\sigma$ is  the image of the diagonal immersion of the group $G_0$   in $G.$ In this case, $L= H \cap K$ is the image of the diagonal immersion of $K_0.$  The associated pair  to $(G,H)$ is $(G, H_0)$ with $H_0$ the image of the immersion of $G_0$ in $G$ via the map $x \rightarrow (x, \theta_0 x).$   Thus, the pair $(H_0, L)$ is isomorphic to the symmetric pair $(G_0, K_0).$    Then, we set $K_1 := \text{diagonal \, subgroup \, of\,} (Z_{K_0} \times Z_{K_0})=:\Delta (Z_{K_0})$ and hence, $K_2= Z_2 ( K_0^s \times K_0^s),$ where  $Z_2 $ is a complement in $Z_{K_0} \times Z_{K_0}$ to the diagonal subgroup $\Delta(Z_{K_0}).$   Hence, $U=H\cap T=\Delta (Z_{K_0}) \Delta (T_0^s)=\Delta (T_0)$  and $L \cap K_2 =\Delta (K_0^s). $ We fix a holomorphic system of positive roots  $\Psi_0$ for $\Phi(\g_0, \t_0).$ Hence $\Psi := \Psi_0 \times \{0\} \cup \{0\} \times \Psi_0$ is a holomorphic system of positive roots for $\Phi(\g, \t).$   Let $(\lambda, \phi)$ denote a Harish-Chandra parameter for $G$ dominant for $\Psi.$ Then, $\pi_\lambda^{G_0}$ and $ \pi_\phi^{G_0} $ are holomorphic irreducible square integrable representations of $G_0$  and since $K_1$ is contained in $H,$ the  outer tensor product  $\pi_{(\lambda, \phi)}^G :=\pi_\lambda^{G_0} \boxtimes \pi_\phi^{G_0} $ is an admissible representation of $H$ \cite{kovinv}.  Because of the decomposition $K= Z_2 \Delta(Z_{K_0}) (K_0^s \times  K_0^s) $  the  lowest $K-$type of $\pi_{(\lambda, \phi)}^G$ restricted to $ \Delta(Z_{K_0}) (K_0^s \times  K_0^s) $ is an irreducible representation and the corresponding Harish-Chandra parameter is $((\lambda_{\z_{\k_0}} + \phi_{\z_{\k_0}})/2 + \rho_n, (\lambda_{\z_{\k_0}} + \phi_{\z_{\k_0}})/2 + \rho_n)+ (\lambda_s, \phi_s),$ where we write $\t^\star \ni \lambda=\lambda_{\z_{\k_0}}+ \lambda_s $ with $\lambda_{\z_{\k_0}} \in \z_{\k_0},  \lambda_s \in \t_0^s.$  The branching law for the restriction of $ \pi_{(\lambda_s,\phi_s)}^{K_0^s \times K_0^s}$ to $\Delta(K_0^s)$ is,   $$ res_{L\cap K_2} (\pi_{(\lambda_s,\phi_s)}^{K_0^s \times K_0^s}) = \sum_{\nu \in Spec_{\Delta(K_0^s)}( \pi_{(\lambda_s,\phi_s)}^{K_0^s \times K_0^s})} m^{K_0^s \times  K_0^s , L \cap K_2} ((\lambda_s,\phi_s), \nu) \,\, \pi_{\nu}^{L \cap K_2}.$$ Following the path of the proof of Theorem 1 we show,
\begin{thm} The multiplicity $m^H ((\lambda,\phi), \mu)$ of $\pi_\mu^H$ in $res_H (\pi_{(\lambda, \phi)}^G) $ is given by
$$  \sum_{\nu \in Spec_{\Delta(K_0^s)}( \pi_{(\lambda_s,\phi_s)}^{K_0^s \times K_0^s} )} m^{G_0, K_0}(((\lambda_{\z_{\k_0}}+\phi_{\z_{\k_0}})/2, \nu), \mu) \, m^{K_0^s \times K_0^s , \Delta(K_0^s)} ((\lambda_s, \phi_s), \nu).$$

\end{thm}

For this case, as in section 4, we may consider   the $L-$invariant  subspaces $\mathcal L_{(\lambda,\phi)}, $ and $ \mathcal U(\h_0)(\mathcal V_{\lambda +\rho_n}^{K_0} \boxtimes V_{\phi +\rho_n}^{K_0}. $  When we assume that both representations are scalar holomorphic,   theorem 3 yields that both subspaces   are equivalent as representations of $L.$  However,  they are not equal, because hypothesis ii) in proposition 4 holds. In fact,
$[[(Y_{-\beta}, Y_{-\beta)}), (Y_{\beta}, -Y_{\beta)})], (Y_{\beta}, -Y_{\beta)})]\not= 0$ for every root $\beta \in \Psi_0.$

\section{Notation for Table I, II and III.}

In this section we complete the notation for the objets presented in the three
 tables and we do a case by case verification of (4.4).

\noindent
$\mathbf{I-3} \,\, (\mathfrak{so}(2m,2n), \mathfrak {so}(2m,2k)+\mathfrak{so}(2m,2n-2k)), m\geq 2, n\geq 2.$ \\ For a suitable orthogonal basis $\epsilon_1, \dots, \epsilon_m, \delta_1, \dots, \delta_n$  of $i\t_\mathbb R^\star $ $$\Delta = \{(\epsilon_i \pm \epsilon_j), 1 \leq i < j\leq m \}\cup \{(\delta_r \pm \delta_s), 1 \leq r < s \leq n \} $$ $$\Phi_n = \{ \pm (\epsilon_r \pm \delta_s), r=1, \dots,m, s=1,\dots,n \}.$$ $\Phi(\l,\t)=\{\pm(\epsilon_i \pm \epsilon_j), 1 \leq i < j\leq m \}\cup \{\pm (\delta_r \pm \delta_s), 1 \leq r < s \leq k \}, $  $\Phi_n(\h, \t)=\{ \pm (\epsilon_r \pm \delta_s), r=1, \dots,m, s=1,\dots,k \}.$  The systems of positive roots $\Psi_\lambda$ so that $\pi_\lambda^G$ is an admissible representation of $H$ are the systems $\Psi_{\pm}$ associated to the lexicographic orders $$ \epsilon_1> \dots> \epsilon_m \,\, > \delta_1> \dots> \delta_n,  \,\,  \epsilon_1> \dots>\epsilon_{m-1}> -\epsilon_m \,\,> \delta_1> \dots>\delta_{n-1}> -\delta_n $$
Here, $\k_1(\Psi_{\pm})=\mathfrak{so}(2m).$

\noindent
$\mathbf{I-4} \,\, (\mathfrak{so}(2m,2n+1), \mathfrak {so}(2m,k)+\mathfrak{so}(2m,2n+1-k)), m\geq 2, n\geq 2.$ For a suitable orthogonal basis $\epsilon_1, \dots, \epsilon_m, \delta_1, \dots, \delta_n$  of $i\t_\mathbb R^\star $   $$\Delta = \{(\epsilon_i \pm \epsilon_j), 1 \leq i < j\leq m \}\cup \{(\delta_r \pm \delta_s), 1 \leq r < s \leq n \} \cup \{\delta_j, j=1,...,n\}$$ $$ \Phi_n = \{ \pm (\epsilon_r \pm \delta_s), r=1, \dots,m, s=1,\dots,n \} \cup \{ \pm \epsilon_j, j=1,\dots,m \}.$$ The systems of positive roots $\Psi_\lambda$ so that $\pi_\lambda^G$ is an admissible representation of $H$ are the systems $\Psi_{\pm}$ associated to the lexicographic orders $$ \epsilon_1> \dots> \epsilon_m> \delta_1> \dots> \delta_n  ,\,  \epsilon_1> \dots>\epsilon_{m-1}> -\epsilon_m> \delta_1> \dots> \delta_{n-1} > -\delta_n. $$
Here, for $M \geq 3$ $\k_1(\Psi_{\pm})=\mathfrak{so}(2m).$ For $m=2, \k_1(\Psi_\pm)=\mathfrak{su}_2(\epsilon_1 \pm \epsilon_2).$

\noindent
$\mathbf{I-5} \,\, (\mathfrak{so}(4,2n), \mathfrak {u}(2,n)_1),  n\geq 3.$
 $\h_\mathbb R=\u(2,n)_{1}$ has for root system $\Phi (\h,\t)=\{\pm(\epsilon_1 -\epsilon_2)\} \cup \{(\delta_k -\delta_s), k\not= s\} \cup \{\pm(\epsilon_1 -\delta_j), \pm (\epsilon_2 -\delta_j), j=1,\dots,n\}.$ $(\h_0)_\mathbb R \sim \u(2,n)$ has for root system $\Phi (\h_0,\t)=\{\pm(\epsilon_1 -\epsilon_2)\} \cup \{(\delta_k -\delta_s), k\not= s\} \cup \{\pm(\epsilon_1 +\delta_j), \pm (\epsilon_2 +\delta_j), j=1,\dots,n\}.$ The system of positive roots $\Psi_\lambda$ so that $\pi_\lambda^G$ is an admissible representation of $H$ is the systems $\Psi_{1,-1}$ associated to the lexicographic order $$ \epsilon_1 >  -\epsilon_2 > \delta_1> \dots> \delta_{n-1} >
  -\delta_n  $$ Here, $\k_1=\mathfrak{su}_2(\epsilon_1 -\epsilon_2).$

\noindent
$\mathbf{I-6} \,\, (\mathfrak{so}(4,2n), \mathfrak {u}(2,n)_2),  n \geq 3.$   $\h_\mathbb R=\u(2,n)_{2}$ has for root system $\Phi (\h,\t)=\{\pm(\epsilon_1 +\epsilon_2)\} \cup \{(\delta_k -\delta_s), k\not= s\} \cup \{\pm(\epsilon_1 +\delta_j), \pm (\epsilon_2 -\delta_j), j=1,\dots,n\}.$ $(\h_0)_\mathbb R \sim \u(2,n)$ has for root system $\Phi (\h_0,\t)=\{\pm(\epsilon_1 +\epsilon_2)\} \cup \{(\delta_k -\delta_s), k\not= s\} \cup \{\pm(\epsilon_1 -\delta_j), \pm (\epsilon_2 +\delta_j), j=1,\dots,n\}.$ The system of positive roots $\Psi_\lambda$ so that $\pi_\lambda^G$ is an admissible representation of $H$ is the systems $\Psi_{11}$ associated to the lexicographic order $$ \epsilon_1> \epsilon_2 > \delta_1> \dots> \delta_n  $$ Here, $\k_1=\mathfrak{su}_2(\epsilon_1 +\epsilon_2).$

\noindent
$\mathbf{I-7} \,\, (\mathfrak{so}(4,4), \mathfrak {u}(2,2)_{xy}).$  $\h_\mathbb R =\u(2,2)_{11}= \u(2,2)_1$ has for root system $\Phi (\h,\t)=\{\pm(\epsilon_1 -\epsilon_2)\} \cup \{\pm(\delta_1 -\delta_2) \} \cup \{\pm(\epsilon_1 -\delta_j), \pm (\epsilon_2 -\delta_j), j=1,2\}.$ $(\h_0)_\mathbb R=w\u(2,2)_{11}$ has for root system $\Phi (\h_0,\t)=\{\pm(\epsilon_1 -\epsilon_2)\} \cup \{\delta_1 -\delta_2, \} \cup \{\pm(\epsilon_1 +\delta_j), \pm (\epsilon_2 +\delta_j), j=1,2\}.$  We set $w_{\epsilon, \delta}= S_{\epsilon_1 -\delta_1}S_{\epsilon_2 -\delta_2}. $  $w_{\epsilon \delta}$ normalizes $\Delta$ and switches  epsilon's in delta's.  The systems of positive roots $\Psi_\lambda$ so that $\pi_\lambda^G$ is an admissible representation of $\u (2,2)_{11}$ are the systems $\Psi_{1, -1} , w_{\epsilon \delta} \Psi_{1,-1}. $ Where, $\Psi_{1, -1}$ is the system associated to the lexicographic order $$ \epsilon_1>  -\epsilon_2> \delta_1>  -\delta_2  $$ Here, $\mathfrak k_1(\Psi_{1,-1})=\mathfrak{su}_2(\epsilon_1 -\epsilon_2), \mathfrak k_1(w_{\epsilon \delta}\Psi_{1,-1})=\mathfrak{su}_2(\delta_1 -\delta_2).$

\noindent
$\mathbf{I-8} \,\, (\mathfrak{so}(4,4),\u(2,2)_{12}).$  $\Phi (\h,\t)=\{\pm(\epsilon_1 -\epsilon_2)\} \cup \{\pm(\delta_1 +\delta_2) \} \cup \{\pm(\epsilon_1 -\delta_1), \pm(\epsilon_1 +\delta_2) , \pm (\epsilon_2 -\delta_1), \pm(\epsilon_2 -\delta_2)\}.$ $(\h_0)_\mathbb R=w\u(2,2)_{12}$ has for root system $\Phi (\h_0,\t)=\{\pm(\epsilon_1 -\epsilon_2)\} \cup \{(\delta_1 +\delta_2) \} \cup \{ \pm(\epsilon_1 +\delta_j), \pm (\epsilon_2 +\delta_j), j=1,2\}$   The systems of positive roots $\Psi_\lambda$ so that $\pi_\lambda^G$ is an admissible representation of $\u (2,2)_{12}$ are the systems $\Psi_{1, -1} , w_{\epsilon \delta} \Psi_{1,1}. $  $\Psi_{1, 1}$ is the system associated to the lexicographic order $$ \epsilon_1> \epsilon_2> \delta_1>  \delta_2  $$ Here, $\mathfrak k_1(\Psi_{1,-1})=\mathfrak{su}_2(\epsilon_1 -\epsilon_2), \mathfrak k_1(w_{\epsilon \delta}\Psi_{1,1})=\mathfrak{su}_2(\delta_1 +\delta_2).$

\noindent
$\mathbf{I-9} \,\, (\mathfrak{so}(4,4),\u(2,2)_{21}=\u(2,2)_2).
$  $\Phi (\h,\t)=\{\pm(\epsilon_1 +\epsilon_2)\} \cup \{\pm(\delta_1 -\delta_2) \} \cup \{\pm(\epsilon_1 -\delta_j), \pm (\epsilon_2 +\delta_j), j=1,2\}.$ $(\h_0)_\mathbb R=w\u(2,2)_{21}$ has for root system $\Phi (\h_0,\t)=\{\pm(\epsilon_1 +\epsilon_2)\} \cup \{(\delta_1 -\delta_2), \} \cup \{ \pm(\epsilon_1 +\delta_j),  \pm (\epsilon_2 -\delta_j), j=1,2)\}.$   The systems of positive roots $\Psi_\lambda$ so that $\pi_\lambda^G$ is an admissible representation of $\u (2,2)_{21}$ are the systems $\Psi_{1, 1} , w_{\epsilon \delta} \Psi_{1,-1}. $ Here, $\k_1(\Psi_{11})=\mathfrak{su}_2(\epsilon_1 +\epsilon_2), \mathfrak k_1(w_{\epsilon \delta}\Psi_{1,-1})=\mathfrak{su}_2(\delta_1 -\delta_2).$

\noindent
$\mathbf{I-10} \,\, (\mathfrak{so}(4,4),\u(2,2)_{22}).$  $\Phi (\h,\t)=\{\pm(\epsilon_1 +\epsilon_2)\} \cup \{\pm(\delta_1 +\delta_2) \} \cup \{\pm(\epsilon_1 -\delta_1), \pm(\epsilon_1 +\delta_2) , \pm (\epsilon_2 +\delta_1), \pm(\epsilon_2 -\delta_2)\}.$ $(\h_0)_\mathbb R=w\u(2,2)_{22}$ has for root system $\Phi (\h_0,\t)=\{\pm(\epsilon_1 +\epsilon_2)\} \cup \{(\delta_1 +\delta_2), \} \cup \{ \pm(\epsilon_1 +\delta_1), \pm(\epsilon_1 -\delta_2) , \pm (\epsilon_2 -\delta_1), \pm(\epsilon_2 +\delta_2)\}$   The systems of positive roots $\Psi_\lambda$ so that $\pi_\lambda^G$ is an admissible representation of $\u (2,2)_{22}$ are the systems $\Psi_{1, 1} , w_{\epsilon \delta} \Psi_{1,1}. $

Note for all of the pairs $(\mathfrak{so}(4,2n), \u (2,n))$ we have $w \h = \h_0$ where $w$ is an automorphism of $\mathfrak{so}(4+2n, \mathbb C)$ which extends the linear  operator of $\t^\star$ so that
 $ w(\epsilon_j)=\epsilon_j, j=1,2, w(\delta_k)=-\delta_k, k=1, \dots,n.$

\noindent
$\mathbf{I-11} \,\,  (\mathfrak{sp}(m,n), \mathfrak {sp}(m,k)+\mathfrak{sp}(n-k)).$  For  some orthogonal basis $\epsilon_1, \dots, \epsilon_m, \delta_1, \dots, \delta_n$  of $i\t_\mathbb R^\star $
\begin{multline*} \Delta = \{(\epsilon_i \pm \epsilon_j), 1 \leq i < j\leq m \} \cup \{  2\epsilon_j, j=1,\dots,m \} \\ \cup \{ (\delta_r \pm \delta_s), 1 \leq r < s \leq n \} \cup \{ 2\delta_j, j=1,\dots,n \}
 \end{multline*}
 $$\Phi_n = \{ \pm (\epsilon_r \pm \delta_s), r=1, \dots,m, s=1,\dots,n \} .$$ The system of positive roots $\Psi_\lambda$ so that $\pi_\lambda^G$ is an admissible representation of $H$ is the system $\Psi_{+}$ associated to the lexicographic order $$ \epsilon_1> \dots> \epsilon_m> \delta_1> \dots> \delta_n. $$
Here, $\k_1(\Psi_+)=\mathfrak{sp}(m).$

\noindent
$\mathbf{I-12},...,\mathbf{I-16}  $   For all of these exceptional groups, we have  $U=T,$ it follows from the tables in \cite{kosh}, \cite{dvc} that the unique system  $\Psi_\lambda $ so that $\k_1 \subset \l$ is for $ \Psi_\lambda$ so that   $\k_1=\mathfrak{su}(\alpha_m )$  and this happens only  for $\Psi_\lambda=   =\Psi_{BS} $ a Borel de Siebenthal system of positive roots.    By explicit calculations we find that both $\Psi_{BS} \cap \Phi(\h, \t) $ and $\Psi_{BS} \cap \Phi(\h_0, \t)$ are again a Borel de Siebenthal system.

\noindent
$\mathbf{II-4}  $   $ (\mathfrak {so}(2m,2n), \mathfrak {so}(2m, 2k+1)\times \mathfrak {so}(2n-(2k+1))).$
We fix an outer automorphism $ \nu$ of $\g$ so that  $\nu (\epsilon_j)=\epsilon_j, j=1, \dots,m;  \nu (\delta_j) =\delta_j, j=1, \dots, n-1, \nu (\delta_n)=-\delta_n $  and such that the fix point subalgebra of $\nu$  is $\mathfrak{so}(2m,2n-1)+\mathfrak{so}(1).$ We set $h_0=0$ and for $k=1, \dots , n-1 $ let $h_k \in i\t_\mathbb R $ be  so that  $ (\epsilon_j -\epsilon_{j+1} )(h_k)=(\epsilon_m -\delta_1)(h_k)=0, (\delta_i -\delta_{i+1})(h_k) = \delta_{ik}, (\delta_{n-1} +\delta_n)(h_k)=0, i=1,\dots,n-1, j=1,\dots,m-1.$ Let $\sigma_k = Ad(exp(\pi i h_k)) \nu.$ Then, $\sigma_k$  is the involution which gives rise  the pair
$ (\mathfrak {so}(2m,2n), \mathfrak {so}(2m, 2n-(2k+1))\times \mathfrak {so}(2k+1)).$\\
For any $\sigma_k, $ obviously,  $\u^\star $ is the subspace spanned by $\epsilon_1, \dots, \epsilon_m, \delta_1, \dots, \delta_{n-1}.$\\
\begin{multline*}\Phi(\h, \u)=\{ \pm (\epsilon_r \pm\epsilon_s), 1 \leq r < s\leq m \}\cup \{ \pm (\delta_i \pm \delta_j), 1 \leq i<j\leq 2(n-k-1)+1, \,\text{or} \\ \pm (\delta_i \pm \delta_j),  2n-2k \leq i<j\leq n-1\} \cup \{ \pm \epsilon_r, r=1, \dots, m\}\cup \{ \pm \delta_i, i=1, \dots, n-1\}\\  \cup \{ \pm (\epsilon_p \pm \delta_q), p=1,\dots,m, q=1, \dots, 2n-2k-1 \}.
\end{multline*}

The two system $\Psi_\lambda$ so that $\pi_\lambda^G$ is an admissible representation for $H$ are: $\Psi_+$ the system associated to the lexicographic orders $\epsilon_1 > \dots \epsilon_m > \delta_1 > \dots >\delta_n $ and  $ \Psi_-:=S_{\epsilon_m -\delta_n}S_{\epsilon_m +\delta_n} \Psi_+.$ Here, $\k_1= \mathfrak{so}(2m).$  It is obvious that $\sigma_k (\Psi_{\pm})=\Psi_{\pm}.$ \\
The corresponding systems $\Psi_{H, \lambda}$ are: $\Psi_+$ the system associated to the lexicographic orders $\epsilon_1 > \dots \epsilon_m > \delta_1 > \dots >\delta_{n-1} $ and  $\Psi_- =S_{\epsilon_m} \Psi_+.$ The systems $\Psi_{H_0, \lambda}$ are: the system associated to the lexicographic orders $\epsilon_1 > \dots \epsilon_m > \delta_1 > \dots >\delta_{n-1} $ and the image of this system by $S_{\epsilon_m}.$

\noindent
$\mathbf{II-6}$ \,\,  $ (\e_{6(2)}, \f_{4(4)}). $ Here, $\mathfrak k= \mathfrak{su}_2(\alpha_m) + \mathfrak{su}_6, \mathfrak l= \mathfrak{su}_2(\alpha_m) + \mathfrak{sp}(3)$ Hence, from \cite{dvc} it follows there is a unique  system of positive roots   such that $\k_1 \subset \l.$ The system is the Borel de Siebenthal   $\Psi_{BS}.$  The  Vogan diagram of $\Psi_{BS}$ is

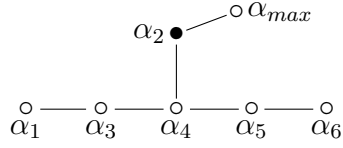
\begin{figure}[h]
\begin{tikzpicture}
\draw (0.05,0) --(0.94,0);  \draw (1.04,0) --(1.97,0); \draw
(2.04,0) -- (2.97,0); \draw (3.03,0) -- (3.96,0) ;
\draw (2,0) -- (2,1);
\draw (2,1) -- (2.8,1.3);
\node[right] at (2.8, 1.3) {$\alpha_{max}$};
\node at (2.8,1.3) {$\circ$};
\node[left] at (2,1) {$\alpha_2$};
\node at (2,1) {$\bullet$};
\node at (2,0) {$\circ$};
\node at (0,0) {$\circ$};
\node at (1,0) {$\circ$};
\node at (3,0) {$\circ$};
\node at (4,0) {$\circ$};
\node[below] at (0,0) {$\alpha_1$};\node[below] at (1,0) {$\alpha_3$};\node[below] at (2,0) {$\alpha_4$};\node[below] at (3,0) {$\alpha_5$}; \node[below] at (4,0) {$\alpha_6$};
\end{tikzpicture}
\end{figure}

\noindent
 Here,   $\alpha_1,  \alpha_3, \dots, \alpha_6$ are the compact simple roots and $\alpha_2$ is the noncompact simple root.   The automorphism  $\sigma$  of $\g$  acts on the simple roots as follows $$\sigma (\alpha_2)=\alpha_2, \,\,  \sigma (\alpha_1)=\alpha_6, \,\,  \sigma( \alpha_3 )=\alpha_5, \,\, \sigma (\alpha_4)=\alpha_4.$$ Hence, $\sigma(\Psi_{BS})=\Psi_{BS}.$  Let $Z_r= \frac{2H_{\alpha_{r}}}{(\alpha_{r}, \alpha_{r})}\in i\t_\mathbb R, r=1,\dots, 6. $ Hence, $\u$ is the subspace spanned by $Z_2, Z_4, Z_1 +Z_6, Z_3 +Z_5. $  and $\t^{-\sigma}$ is spanned by $Z_1 -Z_6, Z_3 -Z_5.$ Let $h_2 \in i\t_\mathbb R^\star $ be so that $\alpha_j(h_2)=\delta_{j 2} $ for $j=1,\dots, 6.$  Hence,   $h_2=\frac{2H_{\alpha_{m}}}{(\alpha_{m}, \alpha_{m})}$ and $\theta= Ad(exp(\pi i h_2)).$ Let $\sigma_2 =\theta \sigma =Ad(exp(\pi i h_2)) \sigma.$ Then, the fix point subalgebra for $\sigma_2$ is isomorphic to $\mathfrak{sp}(1,3)$ and the pair $(\e_{6(2)}, \mathfrak{sp}(1,3))$ is the associated pair to  $ (\e_{6(2)}, \f_{4(4)}). $    The simple roots for $\Psi_{H,\lambda},  \Psi_{H_0, \lambda},$  respectively, are:
\begin{center}
   $\alpha_2, \,\, \alpha_4, \,\,  q_\u(\alpha_3)=q_\u (\alpha_5), \,\,  q_\u(\alpha_1)=q_\u (\alpha_6). $

 $  q_\u(\alpha_2 + \alpha_4 +\alpha_5)=q_\u (\alpha_2 + \alpha_4 +\alpha_3), \,\, q_\u(\alpha_1)=q_\u (\alpha_6), \,\,   q_\u(\alpha_3)=q_\u (\alpha_5), \,\, \alpha_4.   $
\end{center}
The respective Dynkin diagrams are:
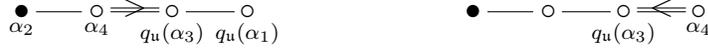
\begin{figure}[h]
\begin{tikzpicture}
 \draw (6.05,0) --(6.94,0);  \draw (7.04,0) --(7.97,0)  ; \draw
(8.04,0) -- (8.97,0) ;  \draw (8.04,0.05) --(8.97,0.05) node[midway]{$<$};
 \node at (8,0) {$\circ$};
\node at (7,0) {$\circ$};
\node at (9,0) {$\circ$};
\node at (6,0) {$\bullet$};
\node[below] at (8,0) {\tiny{$q_\mathfrak u(\alpha_3)$}};\node[below] at (9,0) {\tiny{$\alpha_4$}};
\draw (0.05,0) --(0.94,0);  \draw (1.04,0) --(1.97,0)  ; \draw
(2.04,0) -- (2.97,0);  \draw (1.04,0.05) --(1.97,0.05) node[midway]{$>$};
 \node at (2,0) {$\circ$};
\node at (0,0) {$\bullet$};
\node at (1,0) {$\circ$};
\node at (3,0) {$\circ$};
\node[below] at (0,0) {\tiny{$\alpha_2$}};\node[below] at (1,0) {\tiny{$\alpha_4$}};\node[below] at (2,0) {\tiny{$q_\mathfrak u(\alpha_3)$}};\node[below] at (3,0) {\tiny{$q_{\mathfrak u}(\alpha_1)$}};
\end{tikzpicture}
\end{figure}

   An observation which follows from inspection of the Vogan diagram for each $\Psi_\lambda$ is:  for   a  compact simple root $\alpha$ in $\Psi_\lambda$ we have $\alpha + \sigma (\alpha)$ is not a root.

\smallskip

\noindent
Next, we determine the holomorphic systems which satisfy the hypothesis in i) for proposition 4.

\smallskip

\noindent
$\mathbf{ III-1, III-2}$\,\,$(\mathfrak{su}(m,n), \mathfrak{su}(k,l) +\mathfrak{su}(m-k,n-l)+\u(1)).$ $\h_0=\mathfrak{su}(k,n-l)+\mathfrak{su}(m-k,l)+\u(1).$ We follow  notation in $\mathbf{I-1}.$ For $1 \leq k <m \,\, \text{and} \,\, 1 \leq l< n,$  a nonzero element of $[[\h \cap \p^-, \h_0 \cap \p^+],\h_0 \cap \p^+]$ is
$[[Y_{-\epsilon_1 +\delta_1}, Y_{\epsilon_1 -\delta_{l+1}}], Y_{\epsilon_{k+1} -\delta_1}].$ A nonzero element in $[[\h_0 \cap \p^-, \h \cap \p^+],\h \cap \p^+]$ is  $[[Y_{-\epsilon_1 +\delta_{l+1}}, Y_{\epsilon_1 -\delta_{1}}], Y_{\epsilon_{k+1} -\delta_{l+1}}].$ When  $k=m, 1 \leq l <n, $ we fix  $1\leq i \leq m, 1 \leq j \leq l, 1 \leq a \leq m, l+1 \leq b \leq n, 1 \leq r \leq m, l+1 \leq s \leq n,$ then $[[Y_{-\epsilon_i +\delta_{j}}, Y_{\epsilon_a -\delta_{b}}], Y_{\epsilon_{r} -\delta_{s}}]=[Y_{\delta_j -\delta_{b}}, Y_{\epsilon_{r} -\delta_{s}}]=0.$  Thus,  $[[\h \cap \p^-, \h_0 \cap \p^+],\h_0 \cap \p^+]=\{0\}.$ It readily follows $[[\h_0 \cap \p^-, \h \cap \p^+],\h \cap \p^+]=\{0\}$  as well the analysis when $l=n, k<m.$

\noindent
$\mathbf{ III-3} $ $(\mathfrak{so}(2,2n), \mathfrak{so}(2,2k)+\mathfrak{so}(2(n-k))), n\geq 2.$ We use  notation in $\mathbf{I-3}.$ Here, $[[\h \cap \p^-,\h_0 \cap \p^+ ], \h_0 \cap \p^+] \ni [[Y_{-\epsilon_1 +\delta_{1}},Y_{\epsilon_{1}-\delta_n} ], Y_{\epsilon_1+\delta_n}]\not=0, $ and
$[[\h_0 \cap \p^-,\h \cap \p^+ ],\h \cap \p^+] \ni [[Y_{-\epsilon_1 -\delta_{n}},Y_{\epsilon_{1}-\delta_1} ], Y_{\epsilon_1+\delta_1}]\not=0.$

 \noindent
 $\mathbf{III-4}$ $(\mathfrak{so}(2m,2), \u(m,1)).  $ The holomorphic system $\Psi$ is so that $\Psi \cap \Phi(\l,\t)=\{\epsilon_j -\epsilon_k, 1 \leq j <k\leq m\}, \Psi \cap \Phi_n(\h,\t)=\{\delta_1 -\epsilon_j, 1 \leq j \leq m\}.$ $\Psi \cap \Phi_n(\h_0,\t)=\{\delta_1 +\epsilon_j, 1 \leq j \leq m\}.$ Then,
   $[[Y_{-(  \delta_1-\epsilon_r)}, Y_{\epsilon_k +\delta_1}],Y_{\epsilon_a +\delta_1}]=0,$ and $[[Y_{-(  \delta_1+\epsilon_r)}, Y_{ \delta_1-\epsilon_k }],Y_{-\epsilon_a +\delta_1}]=0.$  Hence,$ [[\h \cap \p^-,\h_0 \cap \p^+ ], \h_0 \cap \p^+]= [[\h_0 \cap \p^-,\h \cap \p^+ ],\h \cap \p^+]=\{0\}. $

\noindent
$\mathbf{ III-5} $ $(\mathfrak{so}(2,2n+1), \mathfrak{so}(2,k)+\mathfrak{so}(2n+1-k)).$  Notation as  in $\mathbf{I-4}.$
For  $ k=1,$  $ [[\h \cap \p^-,\h_0 \cap \p^+ ], \h_0 \cap \p^+] \ni [[Y_{-\epsilon_1 },Y_{\epsilon_{1}+\delta_n} ], Y_{\epsilon_1-\delta_n}]\not=0, $ and $[[\h_0 \cap \p^-,\h \cap \p^+ ],\h \cap \p^+] \ni [[Y_{-\epsilon_1 +\delta_{n}},Y_{\epsilon_{1}} ], Y_{\epsilon_1}]\not=0 .$ $k=2n$ is symmetric to $k=1.$ For $2 \leq k \leq 2n-1,$  $[[\h \cap \p^-,\h_0 \cap \p^+ ], \h_0 \cap \p^+] \ni [[Y_{-\epsilon_1 +\delta_{1}},Y_{\epsilon_{1}-\delta_n} ], Y_{\epsilon_1+\delta_n}]\not=0, $ and $ [[\h_0 \cap \p^-,\h \cap \p^+ ],\h \cap \p^+] \ni [[Y_{-\epsilon_1 +\delta_{n}},Y_{\epsilon_{1}+\delta_1} ], Y_{\epsilon_1-\delta_1}]\not=0. $

\noindent
$\mathbf{ III-6} $\,\, $(\mathfrak{so}^\star(2n), \mathfrak{u}(m,n-m))$ $1 \leq m <n, n\geq 3,$ $\h_0=\mathfrak{so}^\star(2m) \times \mathfrak{so}^\star(2(n-m))).$ The holomorphic system $\Psi$ is so that $\Psi \cap \Phi(\l,\t)=\{ (\epsilon_i -\epsilon_j), 1\leq i < j \leq m \} \cup  \{ (\epsilon_i -\epsilon_j), m+1\leq i < j \leq n \}.$  $\Psi \cap \Phi_n(\h,\t)=\{ (\epsilon_i +\epsilon_j), 1\leq i \leq m < j \leq n \},$ $\Psi \cap \Phi_n(\h_0,\t)=\{ (\epsilon_i +\epsilon_j), 1\leq i
 < j \leq m \}\cup \{ (\epsilon_i +\epsilon_j), m+1\leq i < j \leq n \}.$ For $m=1, 2 \leq j,r,s,a,b $ we have $ [[Y_{-\epsilon_1 -\epsilon_{j}},Y_{\epsilon_{r}+\epsilon_{s}} ], Y_{\epsilon_a+\epsilon_b}] =0,$  hence $ [[\h \cap \p^-,\h_0 \cap \p^+ ], \h_0 \cap \p^+]=\{0\}.$ Also, $ [[Y_{-\epsilon_r -\epsilon_{s}},Y_{\epsilon_{1}+\epsilon_{j}} ], Y_{\epsilon_1+\epsilon_b}] =0,$ whence $ [[\h_0 \cap \p^-,\h \cap \p^+ ],\h \cap \p^+]=\{0\}.$  Analogously for $m=n-1$ we obtain $ [[\h \cap \p^-,\h_0 \cap \p^+ ], \h_0 \cap \p^+]=[[\h_0 \cap \p^-,\h \cap \p^+ ],\h \cap \p^+]=\{0\}.$ For $2\leq m \leq n-2, n\geq 4,$ we have $ [[\h \cap \p^-,\h_0 \cap \p^+ ], \h_0 \cap \p^+] \ni  [[Y_{-\epsilon_1 -\epsilon_{n}},Y_{\epsilon_{1}+\epsilon_{2}} ], Y_{\epsilon_{n-1}+\epsilon_n}] \not= 0,$  and $ [[\h_0 \cap \p^-,\h \cap \p^+ ],\h \cap \p^+] \ni [[Y_{-\epsilon_1 -\epsilon_{2}},Y_{\epsilon_{1}+\epsilon_{n}} ], Y_{\epsilon_2 + \epsilon_{n-1} }] \not= 0.$

   \noindent
$\mathbf{ III-7} $\,\, $(\mathfrak{sp}(n, \mathbb R), \mathfrak{u}(m,n-m))$ $1 \leq m <n, n\geq 3,$  $\h_0= \mathfrak{sp}(m,\mathbb R) \times \mathfrak{sp}(n-m, \mathbb R)).$  The holomorphic system $\Psi$ is so that
$\Psi \cap \Phi(\l,\t)=\{ (\epsilon_i -\epsilon_j), 1\leq i < j \leq m \} \cup  \{ (\epsilon_i -\epsilon_j), m+1\leq i < j \leq n \} ,$   $\Psi \cap \Phi_n(\h,\t)=\{ (\epsilon_i +\epsilon_j), 1\leq i \leq m < j \leq n \},$ and  $ \Psi \cap \Phi_n(\h_0, \t)= \{ (\epsilon_i +\epsilon_j), 1\leq i \leq j \leq m \}
   \cup \{ (\epsilon_i +\epsilon_j), m+1\leq i \leq j \leq n \}. $ For $1 \leq m <n$ we have,
   $ [[\h \cap \p^-,\h_0 \cap \p^+ ], \h_0 \cap \p^+] \ni [[Y_{-\epsilon_1 -\epsilon_{m+1}},Y_{2\epsilon_{m+1}} ], Y_{2\epsilon_1}]\not=0 . $ For $ 1 \leq m \leq n-2,$ we have $ [[\h_0 \cap \p^-,\h \cap \p^+ ],\h \cap \p^+] \ni
   [[Y_{-2\epsilon_1},Y_{\epsilon_1 + \epsilon_{n-1} } ], Y_{\epsilon_1 +\epsilon_n}]\not=0 .$  For $m=n-1, n\geq 2$
$ [[\h_0 \cap \p^-,\h \cap \p^+ ],\h \cap \p^+] \ni [[Y_{-2\epsilon_n},Y_{\epsilon_1 + \epsilon_{n} } ], Y_{\epsilon_{n-1} +\epsilon_n}]\not=0 .$

\noindent
$\mathbf{ III-8} $  $ (\mathfrak e_{6(-14)}, \mathfrak{so}(2,8)+\mathfrak{so}(2)).$ $\beta_\h \in \h, \beta, \gamma \in \h_0, $ if  $[[Y_{-\beta_\h} , Y_\beta ], Y_\gamma]$ were nonzero, then the  coefficient of $\alpha_1$ in $-\beta_\h +\beta+\gamma$ would be 2, a contradiction. Thus, $ [[\h \cap \p^-,\h_0 \cap \p^+ ], \h_0 \cap \p^+]=\{0\}.$ Analogously we obtain $ [[\h_0 \cap \p^-,\h \cap \p^+ ],\h \cap \p^+]=\{0\}.$

\noindent
$\mathbf{ III-9} $ $(\e_{6(-14)}, \mathfrak{su}(2,4)+\mathfrak{su}(2))$, $\h_0= \mathfrak{su}(2,4) +\mathfrak{su}(2).$  In notation of $\mathbf{III-8}.$ $\Psi_n(\h)=\{ \sum d_j \alpha_j : d_3 \in \{0,2\}, d_6=1 \} , $
 $\Psi_n(\h_0)=\{ \sum d_j \alpha_j : d_3= d_6=1 \},$ $\Psi(\l, \t)=\{ \sum d_j \alpha_j : d_3 \in \{0,2\}, d_6=0 \}.$
$\beta=\alpha_1 +\alpha_2+\alpha_3 +\alpha_4 +\alpha_5 +\alpha_6, \beta_{\h_0}= \alpha_3 +\alpha_4 +\alpha_5 +\alpha_6$  we obtain  $ [[\h \cap \p^-,\h_0 \cap \p^+ ], \h_0 \cap \p^+] \ni [[Y_{-\alpha_6},  Y_{\beta_{\h_0}}], Y_{\beta}] \not=0.$  Next, we set $\beta_{\h_0}= \alpha_1 +\alpha_3 +\alpha_4 +\alpha_5 +\alpha_6. $ Then $\beta_{\h}=\alpha_1 +\alpha_2 +2(\alpha_3 +\alpha_4 +\alpha_5)+\alpha_6 = \alpha_3+(\alpha_1+\alpha_2 +\alpha_3 +2(\alpha_4 +\alpha_5)+\alpha_6$ is a root. Hence, $ [[\h_0 \cap \p^-,\h \cap \p^+ ],\h \cap \p^+] \ni [[Y_{-\beta_{\h_0}},  Y_{\beta}], Y_{\alpha_6}] \not=0.$

\noindent
$\mathbf{ III-10} $ $(\e_{6(-14)}, \mathfrak{so}^\star (10)+\mathfrak{so}(2))$, $\h_0= \mathfrak{su}(5,1) +\mathfrak{sl}_2(\mathbb R).$  In notation of $\mathbf{III-8}.$ $\Psi_n(\h)=\{ \sum d_j \alpha_j : d_2=d_6=1 \},  $
 $\Psi_n(\h_0)=\{ \sum d_j \alpha_j : d_2=0, d_6=1 \} \cup \{\alpha_m \}, $ $\Psi(\l, \t)=\{ \sum d_j \alpha_j : d_2=d_6=0 \}, $ we define $\beta_\h=\alpha_1 +\alpha_3 +\alpha_4 +\alpha_5 +\alpha_6, \beta_{\h_0}= \alpha_1 +\alpha_2 +\alpha_3 +\alpha_4 +\alpha_5 +\alpha_6$ since $\alpha_m -\alpha_2 \in \Psi_n(\h_0) $ we obtain  $ [[\h \cap \p^-,\h_0 \cap \p^+ ], \h_0 \cap \p^+] \ni  [[Y_{-\beta_{\h}},  Y_{\beta_{\h_0}}], Y_{\alpha_m -\alpha_2}] \not=0.$   Next,  we set $\beta_{\h_0}= \alpha_1 +\alpha_2 +\alpha_3 +\alpha_4 +\alpha_5 +\alpha_6.$ Then $\alpha_m -\beta_{\h_0}=\alpha_2 +\alpha_3 +2\alpha_4 +\alpha_5 = \alpha_4+(\alpha_2 +\alpha_3 +\alpha_4 +\alpha_5)$ is a root. Hence, $  [[\h_0 \cap \p^-,\h \cap \p^+ ],\h \cap \p^+] \ni [[Y_{-\beta_{\h_0}},  Y_{\alpha_m}], Y_{\alpha_6}] \not=0.$

\noindent
$\mathbf{ III-11} $   $(\e_{7(-25)},\mathfrak{so}^\star(12) +\mathfrak{su}(2)),$ $\h_0= \mathfrak{su}(6,2).$ We use the notation set up by Bourbaki, thus, for the holomorphic system $\alpha_1, \dots \alpha_6$ are the compact simple roots and $\alpha_7$ is the noncompact simple root. $\alpha_1$ is adjacent to the opposite of the maximal root. $\Phi(\l, \t)=\{ \sum d_j \alpha_j : d_2 \in \{0,2\}\,  d_7=0\},$     $\Psi_n(\h,\t)=\{ \sum d_j \alpha_j : d_2=d_7=1\}, $  $\Phi_n(\h_0, \t)=\{ \sum d_j \alpha_j : d_2 \in \{0,2\},d_7=1\}.  $ Let  $\beta_\h=\alpha_3 +(\alpha_1 +\alpha_2 +\alpha_3 +2(\alpha_4 +\alpha_5 +\alpha_6)+\alpha_7)= \alpha_1 +\alpha_2 +2(\alpha_3 +\alpha_4 +\alpha_5 +\alpha_6)+\alpha_7,$ then $-\beta_\h +\alpha_m +\alpha_6+\alpha_7= (\alpha_1 +\alpha_2 +\alpha_3 +2\alpha_4 +\alpha_5) +(\alpha_6+\alpha_7) $  is a root, hence $[[Y_{-\beta_{\h}},  Y_{\alpha_m}], Y_{\alpha_6 +\alpha_7}] \not=0. $     Next, let $\beta_\h= \alpha_1 +\alpha_2 +2(\alpha_3 +\alpha_4 +\alpha_5 +\alpha_6)+\alpha_7, \beta =\alpha_1 +(\alpha_2 +\alpha_3 +2(\alpha_4 +\alpha_5) +\alpha_6+\alpha_7) \in \Psi_n(\h) $ and $-\alpha_m +\beta_\h$ as well as $-\alpha_m +\beta_\h +\beta$ are roots, hence $[[Y_{-\alpha_m },  Y_{\beta_{\h}}], Y_{\beta}] \not=0.$

\noindent
$\mathbf{ III-12} $  $(\e_{7(-25)},\e_{6(-14)} +\mathfrak{so}(2))$ $\h_0= \mathfrak{so}(2,10)+\mathfrak{sl}_2(\mathbb R).$ We use the notation $\mathbf{III-11},$   $\Phi(\l, \t)=\{ \sum d_j \alpha_j : d_1=d_7=0\},$     $\Psi_n(\h,\t)=\{ \sum d_j \alpha_j : d_1=d_7=1\}, $  $\Phi_n(\h_0, \t)=\{ \sum d_j \alpha_j : d_1=0,d_7=1\}\cup \{\alpha_m \}.$  A classical algebras computation yields $ \beta_{\h}=\alpha_1 +(\alpha_2 +\alpha_3 +2\alpha_4 +2 \alpha_5 +\alpha_6 +\alpha_7) \in \Psi_n(\h),$ $ \beta=\alpha_5 +\alpha_6 +\alpha_7  \in \Psi_n(\h_0)$  and that $-\beta_\h +\alpha_m $ as well as  $(-\beta_\h +\alpha_m)+\beta_{\h_0} = \alpha_1 +\alpha_2  +2(\alpha_3 +\alpha_4+\alpha_5+\alpha_6) +\alpha_7 $ are roots. Thus,  $[[Y_{-\beta_{\h}},  Y_{\beta_{\h_0}}], Y_\beta] \not=0.$ Now, $\alpha_m -\alpha_1 \in \Psi_n(\h)$ and   for $\beta_{\h_0}= \alpha_2  +\alpha_3 +\alpha_4+\alpha_5+\alpha_6 +\alpha_7, \beta_\h=\alpha_1+\alpha_2  +\alpha_3 +\alpha_4+\alpha_5+\alpha_6 +\alpha_7$ we obtain $[[Y_{-\beta_{\h_0}}, Y_{\beta_{\h}}], Y_{\alpha_m -\alpha_1}]\not=0.$

\noindent
$\mathbf{ III-13} $ \,\,$(\mathfrak{su}(n,n), \mathfrak{so}^\star(2n)), n \geq 2.$  $\h_0=\mathfrak{sp}(n, \mathbb R)).$ We make use of the notation in $\mathbf{I-1}.$ The outer automorphism $\sigma$ acts on $\t^\star$ as follows $\sigma (\epsilon_j)=-\delta_{n-j+1}, j=1,\dots,n.$ Hence $\u^\star$ is the subspace spanned by $\epsilon_j -\delta_{n-j+1}, j=1,\dots,n$ we note that $\sigma(\beta) \not= \beta$ for at least one noncompact root.

\noindent
$\mathbf{ III-14} $\,\, $(\mathfrak{so}(2,2n), \mathfrak{so}(2,2k+1) +\mathfrak{so}(2n-2k-1)).$  We make use of the notation in $\mathbf{I-3}$ The outer automorphism $\sigma$ acts in $\t^\star$ as $\sigma (\epsilon_1)=\epsilon_1, \sigma(\delta_j)=\delta_j, j=1, \dots, n-1, \sigma(\delta_n)=-\delta_n.$ In this case $\sigma(\beta)\not=\beta$ for at least one noncompact root $ \beta.$

Let $\beta $ be a noncompact root so that $\sigma(\beta)\not= \beta$ then $\h_{q_\u(\beta)}$ (resp. $\h_{q_\u(\beta)}$) is spanned by $X_\beta :=Y_\beta +\sigma(Y_\beta)$ (resp. $V_\beta :=Y_\beta - \sigma(Y_\beta)$). Since, $ [[X_{-\beta},V_\beta], V_\beta] \not=0 $ we have that the hypothesis in condition ii holds for the pairs $\mathbf{III-13, III-14}.$

 \section{Acknowledgements} This note grew up from conversations with Michel Duflo on the paper of Birgit Speh \cite{Speh}.  Part of the work was done during the meeting "Branching problems for unitary representations" organized by The American Institute of Mathematics at Max Plank Institute (Bonn). The author is grateful to the organizers of the workshop, T.  Kobayashi, B. Orsted, B. Speh for their kind invitation to attend to such a nice meeting.

 \end{document}